\documentclass{article}
\usepackage{graphicx} 

\usepackage{amsmath}
\usepackage{amssymb}
\usepackage{xcolor}
\usepackage{graphicx}
\usepackage{authblk}
\usepackage{caption}
\usepackage{subcaption}
\usepackage{appendix}

\usepackage[lmargin=1.0in,rmargin=1.in,tmargin=1.in,bmargin=1.in]{geometry}

\newcommand{\sfrac}[2]{\mathchoice
  {\kern0em\raise.5ex\hbox{\the\scriptfont0 #1}\kern-.15em/
   \kern-.15em\lower.25ex\hbox{\the\scriptfont0 #2}}
  {\kern0em\raise.5ex\hbox{\the\scriptfont0 #1}\kern-.15em/
   \kern-.15em\lower.25ex\hbox{\the\scriptfont0 #2}}
  {\kern0em\raise.5ex\hbox{\the\scriptscriptfont0 #1}\kern-.2em/
   \kern-.15em\lower.25ex\hbox{\the\scriptscriptfont0 #2}}
  {#1\!/#2}}

\newcommand{\half}{{\sfrac{1}{2}}}
\newcommand{\stochF}{{\widetilde{F}}}
\newcommand{\detF}{{\overline{F}}}

\newcommand{\commentout}[1]{}

\usepackage{soul} 

\newcommand{\NewOld}[2]{{{{#1}}{{}}}}
\newcommand{\Add}[1]{{#1}}
\newcommand{\RevOut}[1]{{}}

\usepackage{comment}

\title{A Hybrid Algorithm for Systems of Non-interacting Particles  \Add{with an external potential}}

\author[1]{ Ana Djurdjevac}
\author[2]{Ann Almgren}
\author[2]{John Bell}

\affil[1]{Freie Universität Berlin, Germany}
\affil[2]{Lawrence Berkeley National Laboratory, Berkeley, California, 94720, USA}

\date{}

\begin{document}

\maketitle

\begin{abstract}
    Our focus is on simulating the dynamics of non-interacting particles \Add{including the effects of an external potential}, 
    which, under certain assumptions, can be formally described by the Dean-Kawasaki equation. The Dean-Kawasaki equation can be solved numerically using standard finite volume methods. However, the numerical approximation implicitly requires a sufficiently large number of particles to ensure the positivity of the solution and accurate approximation of the stochastic flux. To address this challenge, we extend hybrid algorithms for particle systems to scenarios where the density is low. The aim is to create a hybrid algorithm that switches from a finite volume discretization to a particle-based method when the particle density falls below a certain threshold. We develop criteria for determining this threshold by comparing higher-order statistics obtained from the finite volume method with particle simulations. We then demonstrate the use of the resulting criteria for dynamic adaptation in both two- and three-dimensional spatial settings \Add{ in the absence of an external potential.  Finally we consider the dynamics when an external potential is included}.
\end{abstract}

\section{Introduction}

The dynamics of a system of non-interacting random-walker particles \Add{under the  influence of an external potential} 
can, in a suitable mathematical setting, be formally described by the Dean-Kawasaki equation.

\Add{Specifically} the microscopic model that we are considering is a random walk of $N$ independent indistinguishable Brownian particles $\{B_t^i\}_{i=1}^N$ \Add{
under the influence of an external potential $V(X)$} 
on a $d$-dimensional torus $\mathbb{T}^d, d \in \{1,2,3\}$ or with Dirichlet or homogeneous Neumann boundary conditions on a rectilinear domain. 
\Add{We assume that the particle positions, $X_i(t)$, evolve according to a Langevin equation
\begin{equation}
    dX_i = - \frac{1}{\gamma} \nabla V(X_i)dt + dB_t^i, \quad i = 1,\dots N.
    \label{eq:langevin}
\end{equation}
}
\Add{We} let $\mu$ be the empirical density \RevOut{that is} defined by
\begin{equation}\label{emp_meas_def}
    \Add{\mu}(x,t) := \frac{1}{N} \sum_{i=1}^N \delta_{B_t^i}(x),
\end{equation}
where $\delta$ is a Dirac distribution. 
Utilizing the It\^{o}'s formula one can formally derive the equation for $\mu$, see \cite{dean_langevin_1996,kawasaki1994stochastic,bressloff2024generalized}. The formal form of the equation is the so-called Dean-Kawasaki equation. In the case of independent particles \Add{including an external potential}, this is a  stochastic partial differential equation \Add{(SPDE)} of the form 
\begin{equation}\label{DK_eq}
    d \mu = \frac{1}{2}\nabla^2 \mu dt + \Add{\frac{1}{\gamma}\nabla \cdot (\mu \nabla V )dt} + \frac{1}{\sqrt{N}} \nabla \cdot \left( \sqrt{\mu} d W\right),
      \end{equation}
where $dW$ is a vector-valued space-time white noise. Note that in the case of interacting particles, the interaction term also appears in the equation.  It is well-known  \cite{konarovskyi2019dean} that due to divergence operator, square root and irregularity of the noise, this equation is highly singular, \Add{ meaning that the standard SPDE theory for singular SPDEs \Add{\cite{hairer2014theory,gubinelli2015paracontrolled}} is not directly applicable}.
Analytically, the only martingale solutions \Add{of Eq. (\ref{DK_eq})} take the form of an atomic measure. Furthermore, such solutions only exist for $N\in \mathbb{N}$. 
\Add{This implies that, from a rigorous mathematical perspective, solutions to Eq. (\ref{DK_eq}) can only be determined by directly modeling the individual particles.}
\RevOut{Nevertheless, standard numerical methods applied to the Dean-Kawasaki equation with truncated noise  
can be used to approximate integral averages of the solution.}
\Add{A common approach to circumvent these theoretical issues is to introduce an high-frequency cutoff for the noise, which then gives a regularized version of  Eq. (\ref{DK_eq})
\begin{equation}\label{eq:DK_with_reg_noise}
        du =\frac{1}{2}\nabla^2 udt+ \Add{\frac{1}{\gamma} \nabla \cdot (u \nabla V)dt} +\frac{1}{\sqrt{N}}\nabla\cdot (\sqrt{u}\; dW^{M}) \;\;\;
\end{equation}
where $dW^M$ represents the regularized noise. We note that in this form, $N$ is no longer required to be an integer.}

\Add{The advantage of using the SPDE approach over discrete particle methods is clear when examining systems with a large number of particles.}
\NewOld{Interest }{In spite of 
theoretical difficulties, interest} in the analysis and simulation of \Add{these types of} regularized Dean-Kawasaki\RevOut{-type} equations has significantly increased over the past decade, with applications ranging from social dynamics \cite{djurdjevac2022feedback,helfmann2021interacting} to fluid dynamics \cite{donev2014reversible}. Equations of the form of Eq. \eqref{DK_eq} also arise in the context of Dynamic Density Functional Theory (DDFT). A detailed discussion on the relationship between deterministic and stochastic DDFT can be found in \cite{archer2004dynamical, donev2014dynamic}.  \RevOut{In particular, the advantage of using the stochastic partial differential equation approach over discrete particle methods is clear when examining systems with a large number of particles.}  

Numerical schemes for Dean–Kawasaki-type equations were considered in \cite{cornalba_dean-kawasaki_2021, cornalba2022regularised,martinez2024finite}.
In \cite{cornalba_dean-kawasaki_2021} the authors show that a system of independent Brownian motions can be approximated in an appropriate metric to arbitrary order by a finite difference or finite element discretization of the Dean-Kawasaki SPDE, \NewOld{when the system contains a sufficiently large number of particles and the initial distribution is bounded away from zero.}{under the assumption that the initial particle distribution must be bounded away from zero.} 
They do not prove positivity of the approximation and allow for the solution becoming negative.  
Their error bound includes a term accounting for the solution becoming negative, which happens rarely  under the assumption $h \gg N^{-1/d}$. They also note that in the \Add{low particle} regime where   $h \leq  N^{-1/d}$\Add{, which corresponds to less than or equal to one particle per mesh cell,} the direct simulation of particles would be less expensive than approximation of the Dean–Kawasaki equation.

In \cite{cornalba2022regularised} 
the authors present a discontinuous Galerkin scheme for the regularized inertial Dean-Kawasaki equation introduced in \cite{cornalba2021well}. In order to address low density regimes and still preserve positivity, they suggest modifying the model by either speeding up the momentum dynamics or adding extra diffusion to the density evolution.
In \cite{martinez2024finite}, a finite element method (FEM) for the Dean-Kawsaki equation is considered. The authors show that a standard Galerkin discretization introduces non-physical correlations in the numerical solution.  They then introduce a  
linear transformation 
to eliminate those artificial correlations in the numerical solution.

The numerical approximation of the Dean-Kawasaki equation using either a finite element or finite volume approach implicitly assumes a sufficiently large number of particles to justify a white noise approximation of the stochastic flux. When the number of particles locally in a simulation is too low this assumption breaks down and the numerical solution can become negative. More subtly, for low local particle density the flux computed using a white noise approximation may not be accurate.  The goal here is to construct a hybrid algorithm that switches to a particle algorithm when the particle density becomes sufficiently small. Note that this scenario naturally arises, for instance, when examining cluster formation, as demonstrated in recent work referenced in \cite{wehlitz2024approximating}.

A hybrid algorithm for particles systems was first proposed by Alexander {\it et al.} \cite{ALEXANDER2002,alexander2003algorithm}.  Their work demonstrates the conservative coupling between \NewOld{an SPDE}{ a stochastic partial differential equation (SPDE)} and a particle algorithm in one dimension.  They show that the resulting hybrid method can capture the mean and variance accurately for both open and close systems. However, their focus is on systems in which the particle density is sufficiently large that the SPDE description is accurate.  Here, we first generalize their approach to multiple spatial dimensions with dynamic adaptation.  We then categorize the behavior of the dynamics of higher-order statistics as the particle density becomes small.  This categorization provides the  criterion for determining where the particle description is required.

The paper is structured as follows. In Section \ref{sect:DK-wp} we briefly recall the well-posedness results for the Dean-Kawasaki equation and its approximation and introduce
the linearized Gaussian approximation. Section \ref{sect:numerics} introduces the finite volume approximation of the Dean-Kawasaki equation 
and the linearized Gaussian approximation. In Section \ref{sec:hybrid} we introduce the hybrid algorithm in one dimension and describe the adaptive hybrid method in higher dimensions.
Section \ref{sect:computations} contains computational examples \Add{that serve to validate the
algorithm in the absence of an external potential and illustrate its capabilities when an external potential is included}.

\section{Dean-Kawasaki equation }\label{sect:DK}

\subsection{Well-posedness of the equation and its approximation}\label{sect:DK-wp}

As mentioned in the introduction, due to singularity of the Dean-Kawasaki equation, well-posedness results are difficult to obtain. In general, some regularization of the noise is needed, since in
\cite{konarovskyi2019dean} the authors showed that the only martingale solution of the Dean-Kawasaki equation is of the form of an empirical measure for $N$ independent Brownian motions given by \eqref{emp_meas_def}. Besides the martingale nature of the solution, they also showed that the equation has no solutions unless $N$ is a positive integer 
and therefore it is potentially not well-suited for numerical analysis or computation. Other solution concepts, such as the stochastic kinetic solution, are also explored. However, since we are interested in computational aspects, this approach will not be addressed here; for more detailed analysis, please refer to \NewOld{\cite{fehrman2024well}}{[1]}. .

In order to address the sensitivity of the solution concept of the original Dean-Kawasaki equation, in \cite{djurdjevac2022weak} the authors suggest a regularized non-linear  SPDE approximation of \eqref{DK_eq}. The modification is constructed to handle situations where the solution approaches zero by replacing the non-Lipschitz square root with a Lipschitz function  \Add{that regularizes it around zero. }
Additionally, the noise is approximated by \NewOld{introducing an }{its} ultraviolet cutoff \NewOld{for }{at} frequencies \NewOld{above }{of} order $M$. 

\NewOld{The authors prove }{It is proved} that, under a suitable coercivity assumption\NewOld{, which relates the regularization of the square root and the noise truncation for a given $N$, }{affecting the choice of parameters $\delta$ and $K$ for a given $N$—} the modified equation is well-posed in the standard weak PDE sense, making it well-suited for numerical computation. Moreover, it can be shown that the solution remains non-negative, conserves mass, and, with an appropriate selection of parameters, error bounds relative to the original particle system can be established.

To assess how the modified equation differs from the  Dean-Kawasaki equation with truncated noise, we note that the modification of the square root  becomes relevant  only when the solution is small, corresponding to a very low-density regime. In the hybrid method, low-density regimes will be handled using particle dynamics. \NewOld{Thus for the hybrid method, we aim to discretize the unmodified SPDE Eq.  (\ref{eq:DK_with_reg_noise}) using a finite volume methods (FVM).}{Therefore, the SPDE that we aim to discretize using the finite volume method (FVM) is the Dean-Kawasaki equation with regularized noise:}

An alternative approximation of the Dean-Kawasaki equation is the so-called linearized Gaussian approximation. The idea is to approximate the fluctuations around the hydrodynamic limit with a Gaussian field. More precisely, if we denote by \Add{$\overline{u}$} 
the hydrodynamic limit when $N$ tends to infinity, then in the case of independent particles, the corresponding PDE is 
\begin{equation}
    \partial_t \overline{{\Add{u}}} = \frac{1}{2}\Add{\nabla^2 } \overline{\Add{u}} \Add{\,+\, \frac{1}{\gamma} \nabla \cdot (\overline{\Add{u}} \nabla V )}\;\;.
    \label{eq:g1}
\end{equation}
The linearized Gaussian approximation of the fluctuations is then given by
\begin{equation}
\Add{du^G} = \frac{1}{2} \nabla^2 {\Add{u^G}} \Add{\,dt}+\Add{\frac{1}{\gamma} \nabla \cdot (\Add{u^G} \nabla V )dt} +\frac{1}{\sqrt{N}} \nabla \cdot (\sqrt{\overline{\Add{u}}} dW^M)  
 \label{eq:g2}
\end{equation}
where $u^G(x,0) = \overline{u}(x,0) = u_0(x)$.  It is well-known 
that the Gaussian approximation 
correctly describes to first-order the central limit fluctuations.
However, the Gaussian approximation fails to predict the large deviation action functional for the particle model, which is accurately predicted by the Dean-Kawaski equation.  Thus, the Gaussian approximation cannot accurately determine the likelihood of observing large deviations from typical hydrodynamic limit behaviour \cite{dirr2020conservative}.

In \cite{cornalba_dean-kawasaki_2021} the authors compared simulations of the Dean-Kawasaki equation and the linearized Gaussian equation. They observed that two models exhibit similar behaviour for the second moments but diverge for the higher moments, which they expected based on the analogy with \cite[Lemma 10]{cornalba_dean-kawasaki_2021}.
As noted above, this work assumed a sufficient number of particles per cell, while here we want to address exactly the case when the number of particles per cell can be arbitrary small locally.

\subsection{Numerics}\label{sect:numerics}

In this section we present a simple finite volume discretization of the Dean-Kawasaki equation in three dimensions. {\Add{For simplicity of the notation, we will present the basic methodology without including an external potential. The modifications needed to include an external potential can be found in Appendix A.}  
The choice of a finite volume approach is motivated by the construction of the hybrid algorithm and is well-suited to the conservation form of the equation.  We also consider a finite volume discretization of the linearized Gaussian approximation given by Eqs. (\ref{eq:g1}-\ref{eq:g2}).  For simplicity, we will consider discretization in time using an Euler-Maruyama scheme; 
investigation of other temporal integrators will be addressed in future work. Without loss of generality, we will assume the volume of the entire domain is one.

Since the focus of the numerics is on the behavior of the system as a function of the number of particles per cell,
we will write the discretization in terms of the number density $q$, defined as the number of particles per unit volume, so that
$\int q(x,t) dx = N$. After re-scaling by setting $q = N \Add{u}$, equation \eqref{eq:DK_with_reg_noise} with truncated noise becomes
\begin{align}
    dq = 
 \frac{1}{2} \nabla^2 q + \nabla \cdot \left(\sqrt{\mathrm{max}(q,0)} \, dW^M \right)
\end{align}
where we introduced $\mathrm{max}(q,0)$ in order to address the issue of taking the square root when the number density becomes negative.
In the context of a finite volume discretization, $M$ corresponds roughly to the number of cells in the overall mesh.

\subsubsection{Finite volume discretization}

We discretize on a domain of size $L_x \times L_y \times L_z$ on a uniform finite volume mesh of size $I\times J \times K$ \RevOut{with $N_c = I \;J \;K$ being  cells} with mesh spacing $\Delta x = \frac{L_x}{I}$, $\Delta y = \frac{L_y}{J}$, $\Delta z = \frac{L_z}{K}$. We then let
\begin{align*}
x_{i} =&( i+\frac{1}{2}) \Delta x \;\;\;\;i= 0,...,I-1 \\
y_{j} =&( j+\frac{1}{2}) \Delta y \;\;\;\;j= 0,...,J-1 \\
z_{k} =&( k+\frac{1}{2}) \Delta z \;\;\;\;k= 0,...,K-1.
\end{align*}
We then define cell  \NewOld{$C_{i,j,k}$}{$(i,j,k)$} to be
\[
C_{i,j,k} = [x_{i-\half},x_{i+\half}] \times  [y_{j-\half},y_{j+\half}]\times  [z_{k-\half},z_{k+\half}]
\]
where $x_{i-\half} = x_i - \frac{\Delta x}{2}$, etc.
The structure of the mesh in 2D is illustrated in Figure \ref{fig:mesh}.
\begin{figure}[h!]
  \centering
  \includegraphics[width=.45\textwidth]{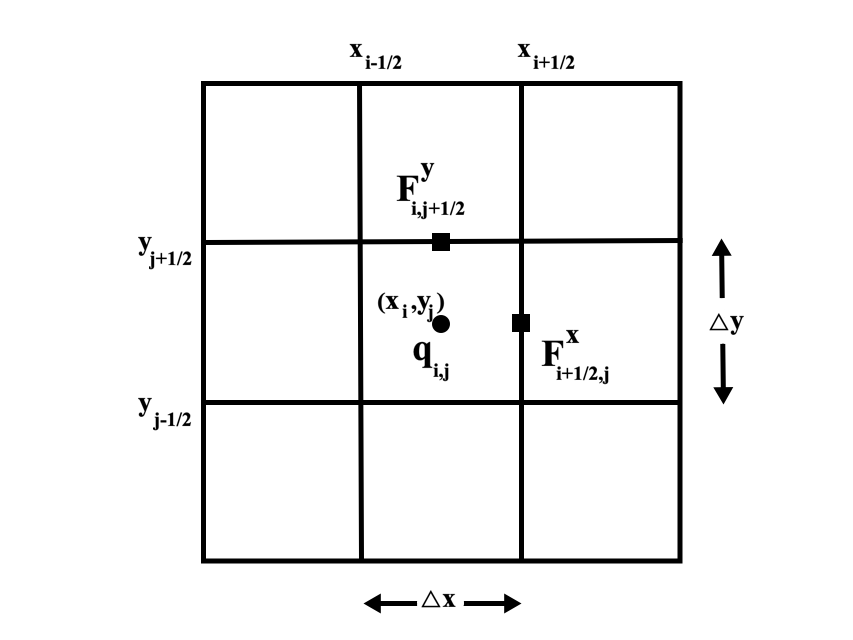}
  \caption{\Add{Sketch of the finite volume mesh in two dimensions. Number density is specified at cell center indicated with a circle. Fluxes are defined
  at edges denoted with squares.
  }  }
  \label{fig:mesh}
\end{figure}

We then let
\begin{equation*}
  q_{i,j,k}(t) \approx  \frac{1}{V_c} \int_{C_{i,j,k}} q(x,t) dx 
\end{equation*}
where $V_c = \Delta x \; \Delta y \; \Delta z$ is the cell 
volume and we  write
\begin{equation*}
    q_{i,j,k}^n \approx q_{i,j,k}(n \Delta t),
\end{equation*}
where $n$ indicates the time step.

We discretize in space and time to obtain 
\begin{equation}\label{eq:cont}
\begin{aligned}
\frac{q^{n+1}_{i,j,k} - q^n_{i,j,k}}{\Delta t} =&
\frac{1}{2} \Big [\frac{q_{i+1,j,k}^n-2q_{i,j,k}^n+ q_{i-1,j,k}^n}{\Delta x^2} +  \frac{q_{i,j+1,k}^n-2q_{i,j,k}^n+ q_{i,j-1,k}^n}{\Delta y^2}  \\ 
+& \frac{q_{i,j,k+1}^n-2q_{i,j,k}^n+ q_{i,j,k-1}^n}{\Delta z^2} \Big ] 
+ \frac{\stochF^x_{i+\half,j,k}-\stochF^x_{i-\half,j,k}}{\Delta x}  \\
+&\frac{\stochF^y_{i,j+\half,k}-\stochF^y_{i,j-\half,k}}{\Delta y}  
+\frac{\stochF^z_{i,j,k+\half}-\stochF^z_{i,j,k-\half}}{\Delta z} 
\end{aligned}
\end{equation}
where the stochastic fluxes are defined as
\[
\stochF_{i+1/2,j,k}^x = \frac{1}{\sqrt{\Delta t V_c}}A(q_{i,j,k},q_{i+1,j,k}) Z_{i+\half,j,k}^x
\]
where \NewOld{$Z_{i+\half,j,k}^x$}{$Z$}
is a Gaussian (normal) distributed random variable with zero mean and unit variance \Add{associated with edge $i+\half,j,k$} and
$A$ is the  averaging operator defined by
\[
A(q_1,q_2) = 
\frac{ \sqrt{\mathrm{max}(q_1,0)} + \sqrt{\mathrm{max}(q_2,0)}}{2}
\]
with
other stochastic fluxes defined analogously.
\Add{This stochastic flux approximates the regularized noise term
$\sqrt{u}\;dW^M$ in Eq. (\ref{eq:DK_with_reg_noise}) where the regularization is given by the mesh spacing.}
We note that we can also define deterministic fluxes
\[
\detF_{i+1/2,j,k}^x = \frac{1}{2} \frac{q_{i+1,j,k}^n - q_{i,j,k}^n}{\Delta x} \;\;\;.
\]
The deterministic and stochastic fluxes can be combined to form a total flux
\begin{equation}
F_{i+\half,j,k}^x = \detF_{i+1/2,j,k}^x + \stochF_{i+1/2,j,k}^x 
\label{eq:totflux}
\end{equation}
with analogous definitions in the other directions.
Using these fluxes we can rewrite Eq. (\ref{eq:cont}) in flux form as
\begin{equation}\label{eq:cont_flux}
\begin{aligned}
\frac{q^{n+1}_{i,j,k} - q^n_{i,j,k}}{\Delta t} =&
\frac{F^x_{i+\half,j,k}-F^x_{i-\half,j,k}}{\Delta x}
+\frac{F^y_{i,j+\half,k}-F^y_{i,j-\half,k}}{\Delta y}   \\ 
+&\frac{F^z_{i,j,k+\half}-F^z_{i,j,k-\half}}{\Delta z} \;\;\; .
\end{aligned}
\end{equation}

\subsubsection*{Linearized Gaussian approximation}

The algorithm described above can easily be modified to treat the linearized Gaussian approximation that corresponds to the equation \eqref{eq:g2}. 
We first define
\begin{equation}\label{eq:cont_gauss}
\begin{aligned}
\frac{\overline{q}^{n+1}_{i,j,k} - \overline{q}^n_{i,j,k}}{\Delta t} =
\frac{1}{2} &\Big [\frac{\overline{q}_{i+1,j,k}^n-2 \overline{q}_{i,j,k}^n+ \overline{q}_{i-1,j,k}^n}{\Delta x^2} +  \frac{\overline{q}_{i,j+1,k}^n-2 \overline{q}_{i,j,k}^n+ \overline{q}_{i,j-1,k}^n}{\Delta y^2}  \\ 
+& \frac{\overline{q}_{i,j,k+1}^n-2\overline{q}_{i,j,k}^n+ \overline{q}_{i,j,k-1}^n}{\Delta z^2} \Big ]\;\;\;.
\end{aligned}
\end{equation}

Using the $\overline{q}$'s, we define Gaussian stochastic fluxes by
\[
\stochF_{i+1/2,j,k}^{G,x} = \frac{1}{\sqrt{\Delta t V_c}}A(\overline{q}_{i,j,k},\overline{q}_{i+1,j,k}) Z_{i+\half,j,k}^x \;\;\; .
\]
We can then define the discretization of the inearized Gaussian approximation as
\begin{equation}\label{eq:cont-Gauss}
\begin{aligned}
\frac{q^{G,n+1}_{i,j,k} - q^{G,n}_{i,j,k}}{\Delta t} =&
\frac{1}{2} \Big [\frac{q_{i+1,j,k}^{G,n}-2q_{i,j,k}^{G,n}+ q_{i-1,j,k}^{G,n}}{\Delta x^2} +  \frac{q_{i,j+1,k}^{G,n}-2q_{i,j,k}^{G,n}+ q_{i,j-1,k}^{G,n}}{\Delta y^2}  \\ 
+& \frac{q_{i,j,k+1}^{G,n}-2q_{i,j,k}^{G,n}+ q_{i,j,k-1}^{G,n}}{\Delta z^2} \Big ] 
+ \frac{\stochF^{G,x}_{i+\half,j,k}-\stochF^{G,x}_{i-\half,j,k}}{\Delta x}  \\
+&\frac{\stochF^{G,y}_{i,j+\half,k}-\stochF^{G,y}_{i,j-\half,k}}{\Delta y}  
+\frac{\stochF^{G,z}_{i,j,k+\half}-\stochF^{G,z}_{i,j,k-\half}}{\Delta z} \;\;\; .
\end{aligned}
\end{equation}

\section{Hybrid discretization}\label{sec:hybrid}

The Dean-Kawasaki \eqref{DK_eq} equation describes the evolution of a large number of independent
Brownian particles \Add{with an external potential}. 
From a computational perspective, when the number of particles in each computational cell is sufficiently large, the stochastic continuum discretization Eq. (\ref{eq:cont})
provides an accurate description of the dynamics.  However, when the number of particles per cell becomes sufficiently small, the approximation is no longer accurate.  Implicit
in Eq. (\ref{eq:cont}) and \eqref{eq:cont-Gauss} is an assumption that the fluctuations in the fluxes can be modeled as white noise. 

The main goal here is to develop a hybrid algorithm that overcomes the restriction of requiring \NewOld{a sufficiently large number of}{ enough} particles  per cell.  The basic idea of the hybrid algorithm is to apply a particle-based method in the regions of the domain where the number density \RevOut{approximation} becomes very small or in the regions where a higher fidelity particle representation is required. The construction of the hybrid algorithm will be based on the adaptive mesh and algorithm refinement (AMAR) approach introduced in \cite{garcia1999adaptive} and extended to stochastic systems in \cite{bell2007algorithm,donev2010hybrid}. Adopting this type of approach, we can treat regions where the particle density is very small (or zero) accurately and ensure positivity of the solution. The algorithm here generalizes the one-dimensional algorithm in \cite{ALEXANDER2002} to multiple spatial dimensions\Add{, includes an external potential,} and incorporates dynamic refinement. 

 AMAR is similar to an adaptive mesh refinement algorithm.  In an adaptive mesh algorithm, regions requiring additional resolution are identified and a finer mesh is used in those regions. In AMAR, one switches to a higher-fidelity model in regions requiring more ``resolution" instead of refining the spatial mesh.
In the AMAR approach, a description of the solution in terms of the SPDE is maintained over the entire domain while particle regions are defined only where needed.  Once the system is initialized the evolution of the system is a three step process.  First, the SPDE region is advanced over the entire domain using Eq. (\ref{eq:cont}).  The SPDE region is then used to supply boundary conditions to advance the particle region. Finally, when both regions have been advanced in time, a synchronization operation is performed to construct a composite solution to the system.

Constructing an 
AMAR hybrid 
involves a number of elements.
We need to define how to map between different representations of the solution.  We need to define a synchronization procedure that corrects results of advancing the SPDE and particle regions sequentially to obtain a composite solution.
Finally, we need to determine where the finer-scale model is required; i.e., when the SPDE model does not provide a sufficiently accurate representation of the system.  One potential criterion here would be simply maintaining positivity of the solution. However, this assumes that the SPDE representation is adequate up until the point when the solution becomes negative.  We will explore this issue in more detail and suggest alternative criteria in the next section.  
The mapping from the particle representation to the SPDE representation is straightforward.  We view the particles as point particles so we can define the SPDE representation as the number of particles in the cell divided by the cell volume.
The mapping from the SPDE representation to a particle representation is not unique.  Here we will view the problem of mapping from SPDE to particles as a conditional sampling problem, namely, we want to generate a random set of particles that are, at least approximately, consistent with the SPDE representation.   Given the SPDE description we can compute the number of particle in each cell.  Those particle are then assigned random locations within the cell where we assume that the particles are uniformly distributed in the cell \Add{ for the case when there is no external potential.  When an external potential is included the sampling needs to be modified, as discussed in Appendix A}.  We note that, in general, the number of particles in a cell may not be an integer. 
If this occurs we use a probabilistic approach to define the particles. Specifically, if a cell contains $\ell + \alpha$ particles where $\ell \, \epsilon \, \mathbb{Z}$ and $\alpha \, \epsilon \, (0,1)$ then we first insert $\ell$ particles as described above. We then choose a random number $\beta$ from the uniform distribution on $[0,1]$. If $\alpha \leq \beta$ then we insert an additional particle; otherwise we do not.  This procedure preserves the expected value of the number of particles in the system but is not strictly conservative. 

\subsection{Hybrid algorithm in one dimension}

We are now ready to discuss the time-step algorithm for the hybrid. For simplicity of the exposition, we will first describe the algorithm in one dimension and assume that the particle region is a contiguous collection of cells from $i_1$ to $i_2$ that does not change with time.  We will then describe the steps needed to extend the overall approach to multiple dimension with dynamic particle patches that consist of different (not necessarily) connected cells. We will also assume that the SPDE and the particle components of the algorithm take the same time step; however, this restriction can also be relaxed.
We assume that the system has been initialized by defining the location of each particle in the particle region and by specifying the number density $q$ in the SPDE region outside the particle region.  
We note that using the mappings defined  above, we can initialize the system from either a pure particle description or from  an SPDE description.  As noted  above, in AMAR the coarse, SPDE representation is maintained throughout the entire domain. Thus, after initializing the particle representation we map the particle solution onto the SPDE solution in the region covered by the particles.

To evolve the system, we first advance the SPDE for a timestep over the entire domain using Eq. (\ref{eq:cont}) and denote the result by $q_i^{n+1,*}$.  To advance the particle region, we need to supply suitable boundary conditions.  The boundary conditions are necessary to capture the effect of particles entering the particle region during a time step.  These boundary conditions are supplied by the SPDE solution.  In the cells adjacent to the particle region, $c_{i_1-1}$ and $c_{i_2+1}$ we compute the number of particles in each  cell from the SPDE solution. We then populate these boundary cells with the specified number of particles with random locations within the cell.  We then advance the particle region by allowing particles in the region to perform a random walk step. \NewOld{Here, we limit the motion of a single particle to be at most $\Delta x$ so that we do not need to account for particles moving more than one cell.  For the time steps used here, $\Delta t$ is sufficiently small that a particle moving more than $\Delta x$ in a time step would be a rare event.  Specifically for the parameters used here (see Sec. \ref{sect:computations}) the probability of a step larger than $\Delta x$ is less than $2 \times 10^{-8}$. Furthermore, empirically we find that when a step length is reduced the reduction is typically less than $0.1\times \Delta x$.  We have found no evidence that these small adjustments in rare cases impact the statistics of the simulations.  (We note that one can enlarge the domain of influence of a particle to allow larger particle steps by enlarging the regions where the SPDE provides boundary data for the particle region.)}{(We note here that we limit the motion of a single particle to be at most $\Delta x$ so that a particle can move at most one cell in a time step. 
We select the time step so that violating this restriction would be extremely rare.)} 
We then define the new particle representation to be the collection of particle in cells $c_{i_1}$ to $c_{i_2}$. (Particles outside the particle region at the end of the step are discarded.)

We now need to synchronize the two representations.  The idea here is that the particle representation is higher fidelity than the SPDE representation so the composite solution should reflect the particle information.  There are two different components to the synchronization.  First, for SPDE cells $i$ that are within the particle region $i_1 \leq i \leq i_2$, we replace the SPDE solution $q_i^{n+1,*}$
with the value obtained from counting the number of particles in the cell \Add{and dividing by the cell volume}, which we denote by $q_i^{n+1}$. 

The other issue we need to address is that the flux used to advance the SPDE cells adjacent to the particle region is inconsistent with the flux computed from the particle region.  To correct this discrepancy, we define
\begin{align}
{q}_{i_1-1}^{n+1} &= q_{i_1-1}^{n+1,*} -\frac{{\Delta t}}{\Delta x} F_{i_1-\half} - \frac{\Delta \mathcal{N}_{i_1-\half}}{\Delta x}  \nonumber \\
&= q_{i_1-1}^{n+1,*}-\frac{\Delta t}{\Delta x}\left[ \frac{q_{i_1}^n-q_{i_1-1}^n}{2 \Delta x} + \widetilde{F}_{i_1-\half} \right ] - \frac{\Delta \mathcal{N}_{i_1-\half}}{\Delta x}  \nonumber \\
 {q}_{i_2+1}^{n+1} &= q_{i_2+1}^{n+1,*}+\frac{{\Delta t}}{\Delta x} F_{i_2+\half} +\frac{\Delta \mathcal{N} _{i_2+1/2}}{\Delta x} \nonumber \\
 &= q_{i_2+1}^{n+1,*}+\frac{\Delta t}{\Delta x}\left [ \frac{q_{i_2+1}^n-q_{i_2}^n}{2 \Delta x} + \widetilde{F}_{i_2+\half} \right ]+\frac{\Delta \mathcal{N} _{i_2+1/2}}{\Delta x}
 \label{eq:reflux}
\end{align}
where $\Delta \mathcal{N}_{i+\half}$ is the net number of particle that crossed edge $i+\half$
from left to right during the time step.  After these adjustments, the SPDE values at
${q}_{i_1-1}^{n+1}$ and ${q}_{i_2+1}^{n+1}$ have been effectively updated using the particle flux instead of the SPDE flux.  For cells away from the particle region 
we simply set $q_i^{n+1}=q_i^{n+1,*}$ to complete the synchronization. The combination of supplying boundary conditions for the particles from the SPDE solution and synchronizing the fluxes using Eq.(\ref{eq:reflux}) effectively creates a composite solution of the system. 
\Add{We note that after the completion of this step, the hybrid integration step is conservative.}
The process is analogous to solving a second-order parabolic
PDE in different subdomains and combining them by matching Dirichlet and Neumann boundary conditions at the interface between the domains.

We note that here we have used the same time step for the SPDE and the particles\RevOut{, where the time step was selected so that a particle moving more than $\Delta x$ during a time step is a rare event}. A smaller time step can be used in the particle region where the particles in the boundary cells for intermediate steps are \NewOld{ populated using the SPDE solution interpolated temporally to the intermediate times}{computed from interpolation of the SPDE solution in time} and the number of particles crossing the boundary is summed over the substeps.

\subsection{Adaptive hybrid in multiple dimensions}

There are a number of additional issues that need to be addressed to extend the basic time step approach discussed above to a dynamically adaptive algorithm in multiple space dimensions.  All of relevant components have analogs in block-structured adaptive mesh refinement, reflecting the overall AMAR approach as a generalization of adaptive mesh refinement (AMR). Some of the details of the algorithm reflect standard practice in AMR.  Further exploiting this relationship, we have implemented the algorithm discussed below using the AMReX framework \cite{AMReX:IJHPCA} for adaptive mesh refinement that provides support for dynamic refinement in parallel and enables the use of GPUs for simulations.  
The basis of block-structured AMR algorithms is a hierarchical representation of the solution at multiple levels of resolution. 
At each level the solution is defined on the union of disjoint data containers at that resolution, each of which represents the solution over a logically rectangular subregion of the domain. 
In the hybrid algorithm both the SPDE region and the particle region are decomposed into a union of disjoint rectilinear grid patches.  On a distributed memory machine, the grid patches are distributed to MPI ranks for parallel execution.

The first issue we consider is the procedure for dynamically modifying the particle region as the solution evolves, which is typically referred to as regridding.  For this discussion, we assume that we have an ``error-estimation" criterion for determining where a particle representation is needed based on the SPDE solution value.  Using the hybrid solution at a given time, we first identify (tag) the cells at the SPDE level where the error criterion indicates that a particle description is needed. As part of this operation we also tag neighbors of cells where the particle representation is needed, to provide a buffer to the particle region.  This reflects the notion that for cells where particles are needed the SPDE flux may not be accurate. We note that this is done over the entire domain including the SPDE representation in the region currently covered by particles. Once we have generated the tagged cells, the cells that are tagged are decomposed into a union of boxes using a clustering algorithm developed by Berger and Rigoutsis \cite{bergerRigoutsos:1991} until the resulting collection of boxes meets a specified efficiency criterion, namely that the fraction of  cells in a box that are tagged is sufficiently large.

Once the grids that define the new particle region have been determined, they need to be filled with data.  Cells that were in the particle region prior to regridding, inherit the particles that were in that cell.  For cells that were not in the particle region, we need to define a particle representation from the SPDE solution.  For this operation, we use the conditional sampling algorithm described above to generate a particle realization that is consistent (approximately) with the SPDE solution.  (We again note that when defining a new particle region from cells with only an SPDE representation, the methodology is not discretely conservative because the SPDE solution may not correspond to an integer number of particles.) The implementation supports regridding at specific intervals specified at input.

We note that essentially the same procedure can be used to initialize the problem. We start by defining the SPDE solution on the entire grid and using the error estimation \Add{criterion} to identify where the particle representation is needed.  The only difference is that at initialization we \Add{have the option to} initialize the particle region directly from a particle description instead of generating random particles based on the SPDE representation.

Although the basic outline of the time step is the same, there are a number of issues that arise because the particle region is now a union of rectilinear  regions in multiple dimensions.  As in the one-dimensional algorithm we first advance the SPDE over the entire domain, including cells covered by the particle domain. 
We note that, in practice, the SPDE region is decomposed into multiple rectilinear subregions.  When implemented in this form, the boundary region for a given grid patch that is in the interior of the SPDE region is filled from the SPDE data on other grids that contain the boundary points.  In addition, at faces at the intersection of two different patches, the random numbers used to compute the stochastic fluxes are synchronized so that the flux at the face is consistent between the two grids.

To advance the particle region we must first provide boundary data for the particle region. 
As in one-dimension, we limit the distance any particle can travel one mesh spacing in each direction and choose the time-step so that a step larger than the mesh spacing would be a rare event.  This implies that to advance the particle solution at a given point, we need  the particle representation in a $3\times 3$ neighborhood  of the point in two dimensions and in a $3 \times 3 \times 3 $ neighborhood in three dimensions.  Thus, we need boundary data for cells $(I,J,K)$ that are not in the particle region but where there is a cell $(i,j,k)$ in the particle region such that
\[
max( | I-i|,|J-j|,|K-k|) \leq 1 \;\;\;.
\]
The necessary boundary cells are filled based on data from the SPDE description.  In particular, from the SPDE solution, we determine the number of particles in the region, using the same probabilistic approach discussed in one dimension.  These particles are then  assigned positions based on sampling from a uniform distribution over the cell, as indicated by the dark blue particles in Figure \ref{fig:Hyb} for the two-dimensional case. \Add{As in one dimension, the sampling need to be modified when there is an exernal potential.}
\begin{figure}[h!]
  \centering
  \includegraphics[width=.45\textwidth]{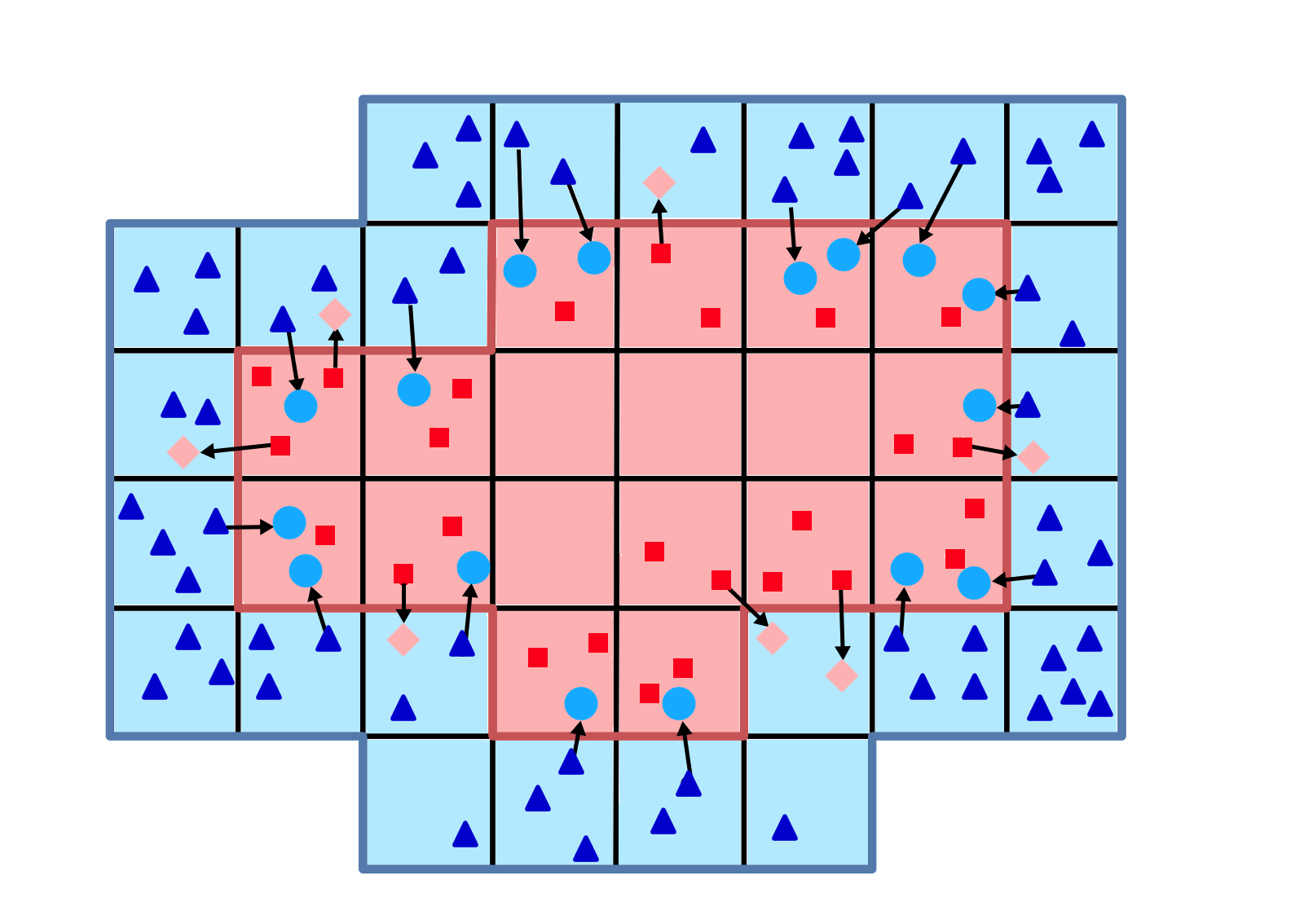}
  \caption{Sketch of hybrid algorithm in two dimensions.  The red-shaded region correspond to the particle region; the blue-shaded region indicates boundary cells needed to advance particle region. {\Add{The dark blue triangles in the blue-shaded region are generated probabilistically to provide a particle configuration consistent with the local SPDE solution. Arrows indicate particles that cross the boundary of the particle region during the time step  that are used to compute particle fluxes into and out of the particle region. Here, the dark blue triangles connected to light blue circles represent particles entering the particle region and the dark red squares connecting to light red diamonds represent particles leaving the particle region.}
  }  }
  \label{fig:Hyb}
\end{figure}
We note that when implemented to operate on rectilinear patches, some of the boundary cells for one patch ($G_1$) may correspond to particle cells in another patch ($G_2$).  In this case the data is copied from $G_2$ to $G_1$ so that the representations are consistent.  Furthermore, some cells can be boundary cells for two different grids.  In this case, we need to ensure that both the number of particles and their locations are the same on both patches.  Finally, we must guarantee that both interior particle cells and boundary particle cells that appear in two different grid patches (including their boundary cells) use the same random numbers to advance the particles. (This is implemented by associating the random numbers as attributes of the particles so that realizations of the same particle in different grid patches are consistent.)

The final step in the hybrid time step is the synchronization of the particle and SPDE
representations. As in one dimension, for SPDE cells that are in the particle region, we replace the SPDE representation by the number density determined by the particle representation; i.e., the number of particles in the cell divided by the cell volume.
The second part of the synchronization is to modify SPDE cells near the particle region so that their updates reflect changes due to the particle update rather than the
SPDE fluxes from SPDE cells covered by the particle region, analogous to Eq.  (\ref{eq:reflux}).  In multiple dimensions with more general particle regions, this correction becomes more complex.  A given cell may be updated with multiple SPDE fluxes. Furthermore, the particle update uses data from a larger set of cells than the SPDE update.  The cells that need to be updated here are cells \NewOld{$C_{I,J,K}$}{ $(I,J,K)$} that are not covered by a particle cell but for which there is as least one cell \NewOld{$C_{i,j,k}$}{ $(i,j,k)$} in the particle region where 
\begin{equation}
max( | I-i|,|J-j|,|K-k|) \leq 1 \;\;\;.
\label{eq:range}
\end{equation}

To describe the algorithm, we first note that the total fluxes from Eq. (\ref{eq:totflux}) can be scaled to give the number of particles transported into the cell across a face during a time step.  For example, for cell \NewOld{$C_{I,J,K}$}{ $(I,J,K)$}
\[
\mathcal{F}^{(I-1,J,K)\rightarrow (I,J,K)} = \Delta t \; \Delta y \; \Delta z \; 
F_{I-\half,J,K}^x 
\]
is the number of particles that transferred from \NewOld{$C_{I-1,J,K}$}{ $(I-1,J,K)$} to \NewOld{$C_{I,J,K}$}{ $(I,J,K)$} and 
\[
\mathcal{F}^{(I+1,J,K) \rightarrow (I,J,K)} = -\Delta t \; \Delta y \; \Delta z \; 
F_{I+\half,J,K}^x 
\]
is the number of particles that transferred from \NewOld{$C_{I+1,J,K}$}{$(I+1,J,K)$} to \NewOld{$C_{I,J,K}$}{$(I,J,K)$}.
We can also define $\Delta N^{(i,j,k) \rightarrow (I,J,K)} $ to be the net number of particles transferred from \NewOld{$C_{i,j,k}$}{ $(i,j,k)$} to \NewOld{$C_{I,J,K}$}{$(I,J,K)$} in the particle update.  The computation of $\Delta N$ is illustrated in Figure \ref{fig:Hyb}.  \NewOld{Arrows that point from dark red squares inside the particle region to light red diamonds particles in the boundary region generate a positive increment to $\Delta N$ for the boundary cell that contains the light red diamonds particle. Arrows that point from dark blue triangles particles to light blue circles inside the particle region decrement $\Delta N$ for the boundary cell containing the dark blue  triangle particle.}
{Arrows that point from red particles inside the particle region to magenta particles in the boundary region generate a positive increment to $\Delta N$ for the boundary cell that contains the magenta particle. Arrows that point from dark blue particles to cyan particles inside the particle region decrement $\Delta N$ for the boundary cell containing the dark blue particle.}
We then define $S_1$ to be the set of particle cells \NewOld{$C_{i,j,k}$}{$(i,j,k)$} that share a face with SPDE cell 
\NewOld{$C_{I,J,K}$}{$(I,J,K)$} and $S_2$ is the set of particle cells that satisfy Eq. (\ref{eq:range}).
We then correct the SPDE solution $q_{I,J,K}^{n+1,*}$ using
\begin{equation*}
q_{I,J,K}^{n+1} = q_{I,J,K}^{n+1,*} - \frac{1}{V_c} \sum_{(i,j,k)\epsilon S^1}
\mathcal{F}^{(i,j,k)\rightarrow (I,J,K)} + \frac{1}{V_c} 
 \sum_{(i,j,k)\epsilon S^2} \Delta N^{(i,j,k) \rightarrow (I,J,K)} \quad _.
\end{equation*}
This completes the synchronization step, giving the final composite solution at the end of the time step.

\section{Computational examples}\label{sect:computations}

In this section, we present several computational examples relating to the behavior of the particle system. We will compare particle simulations \Add{with the } finite volume method, the linearized Gaussian approximation and the hybrid method.  \Add{The initial tests, which consider problems without an external potential, are designed to characterize the properties of the methodology in a simple setting.  These initial simulations are performed on a periodic domain. The last example incorporates an external potential to illustrate the application of the methodology to rare event simulation.}

First, we examine the differences between the particle description and the SPDE description for small particles numbers.  We note that, for the simple Euler-Maruyama discretization presented here, the temporal truncation error leads to an over-prediction of the variance that converges with first-order accuracy in $\Delta t$.  However, skewness and kurtosis are insensitive to the time step.

As noted above, the SPDE description of the particle system is based on an assumption that there are enough particles locally in the system for the stochastic flux of particles to be approximated by white noise.  As the particle number becomes small, this assumption breaks down.  One clear difference is that the PDF of the number of particles per cell from a pure particle algorithm takes on discrete values, namely $k  / V_c$ where $k$ is non-negative integer, while the PDF for the SPDE is continuous.
However, there are additional differences.  
In Figure \ref{fig:pdfs} we show numerical PDFs of systems with 500 and 100 particles at equilibrium 
in one-dimension discretized with 100 \NewOld{cells}{ zones} ($\Delta x = 0.01$) corresponding to an average of 5 and 1 particles per cell, respectively\Add{, using a time step $\Delta t = 3.125\times 10^{-6}$, which is small enough to mitigate the over-prediction of variance discussion above.}  For particle simulations we plot the PDF values as discrete points. For the SPDE we choose bins of size $[ ( k-\half) / V_c , ( k+\half) / V_c]$ which correspond to the width of one particle per cell centered at discrete particle values and plot those values as a continuous distribution.
The PDF from the particle simulations show  asymmetry for the 500 particle case that becomes more significant with 100 particles. 
For both cases, the linearized Gaussian algorithm results in a Gaussian distribution that includes negative values, which becomes more significant \Add{as the number of particles decreases}.  The finite volume algorithm captures some of the asymmetry of the particle algorithm, but the agreement is fairly poor for 100 particles.  In addition, although not as severe as the linearized Gaussian algorithm, the finite volume algorithm also generates negative values for both 5 particles per cell and 1 particle per cell. 
\begin{figure}[h!]
  \centering
 \includegraphics[width=.45\textwidth]{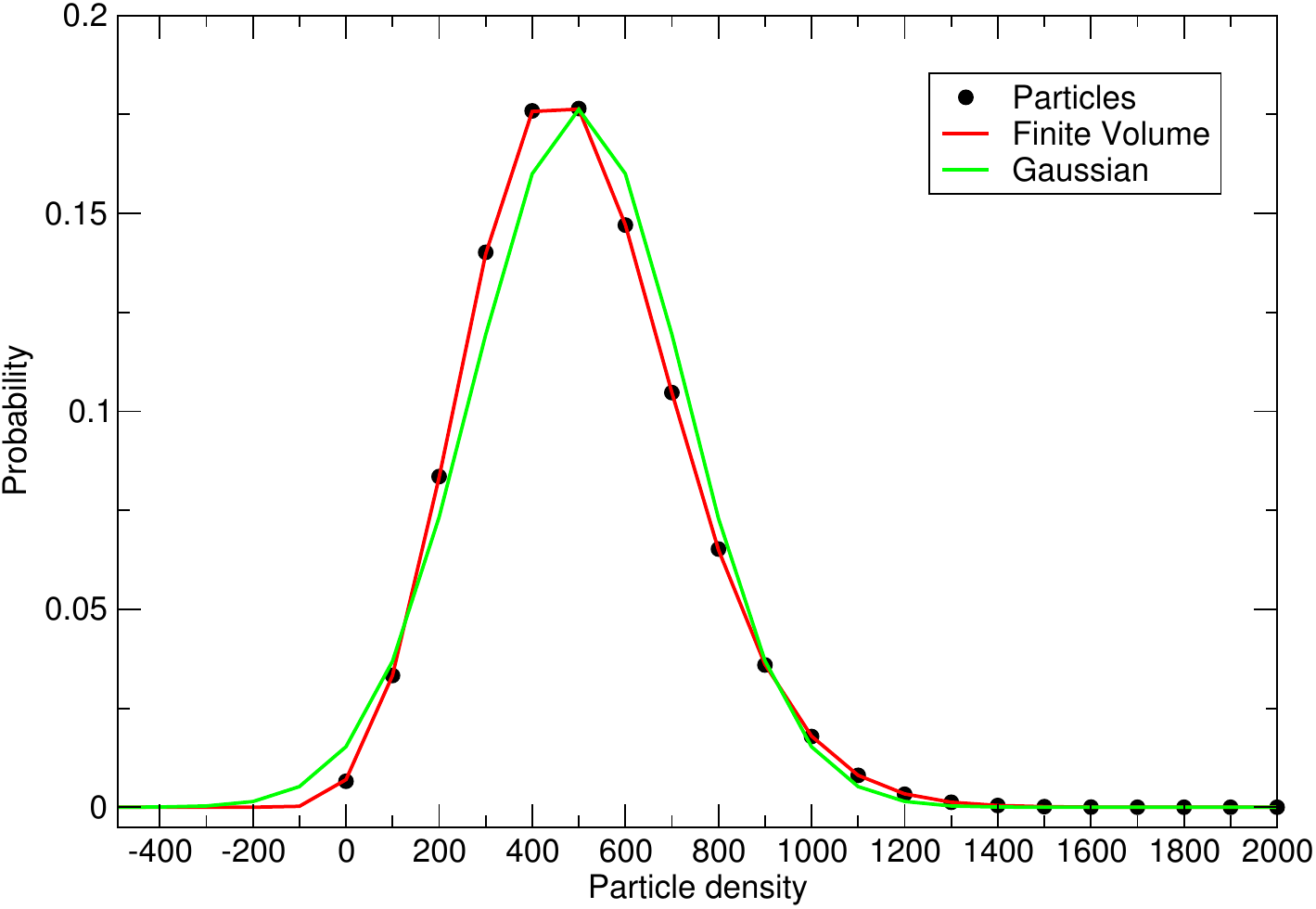}
   \includegraphics[width=.45\textwidth]{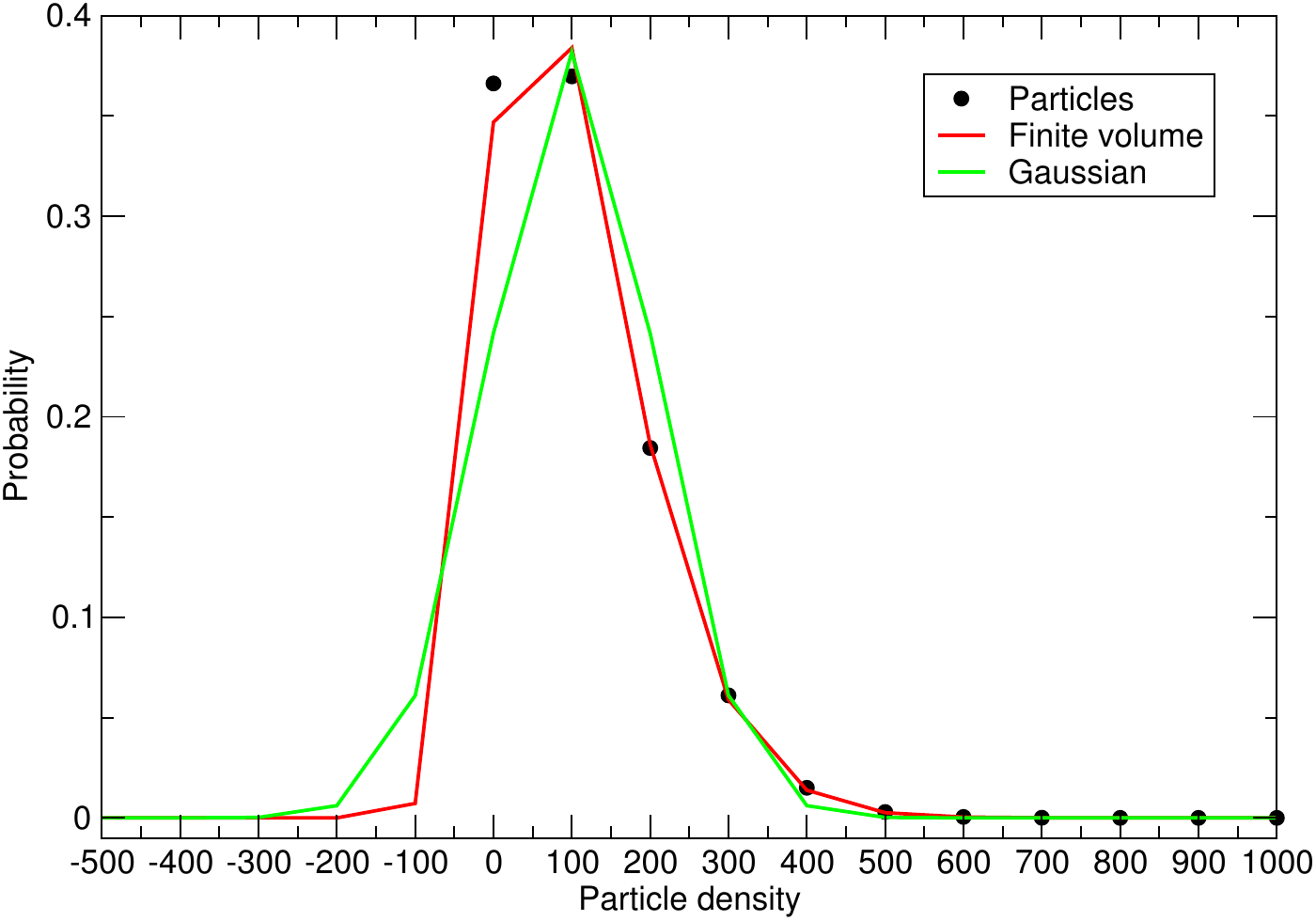}
  \caption{PDFs of particle, finite volume and Gaussian numerical methods at low number densities.  Left is 5 particles per cells; right is 1 particle per cell. \Add{With five points per cell 1.3\% of the cells are negative for Gaussian algorithm and 0.2\% are negative for the finite volume algorithm.
  With one point per cell 15.9\% of the cells are negative for Gaussian algorithm and 12.3\% are negative for the finite volume algorithm.}}
  \label{fig:pdfs}
\end{figure}

To quantify the differences in these distributions, we examine higher-order moments of the distribution.
In Figure \ref{fig:kurt_skew}, we compute skewness and kurtosis of the equilibrium distributions for 
each method as a function of the number of particles per cell, again using a one dimensional problem with 100 \NewOld{cells}{ zones}.  With 50 particles per \NewOld{cell}{ zone}, all of the methods predict a skewness close to zero and a kurtosis of about three, corresponding to a Gaussian distribution.  The Gaussian algorithm retains the Gaussian character (by construction) for all values of particles per cell.  For the particle algorithm and the finite volume method, deviations from Gaussian begin to appear at 20 particles per cell and become significant at 10 particles per cell. \Add{These deviations reflect the underlying multinomial distribution for a finite size domain. For comparison purposes in Figure \ref{fig:kurt_skew} we include the skewness and kurtosis for the Poisson distribution, which approximates the multinomial distribution for large numbers of cells.}   The particle algorithm and the finite volume algorithm agree reasonably well at 10 particles per cell but begin to show \Add{minor} discrepancies at 5 particles per cell that become more pronounced at 1 particle per cell.  
\begin{figure}[h!]
  \centering
  \includegraphics[width=.85\textwidth]{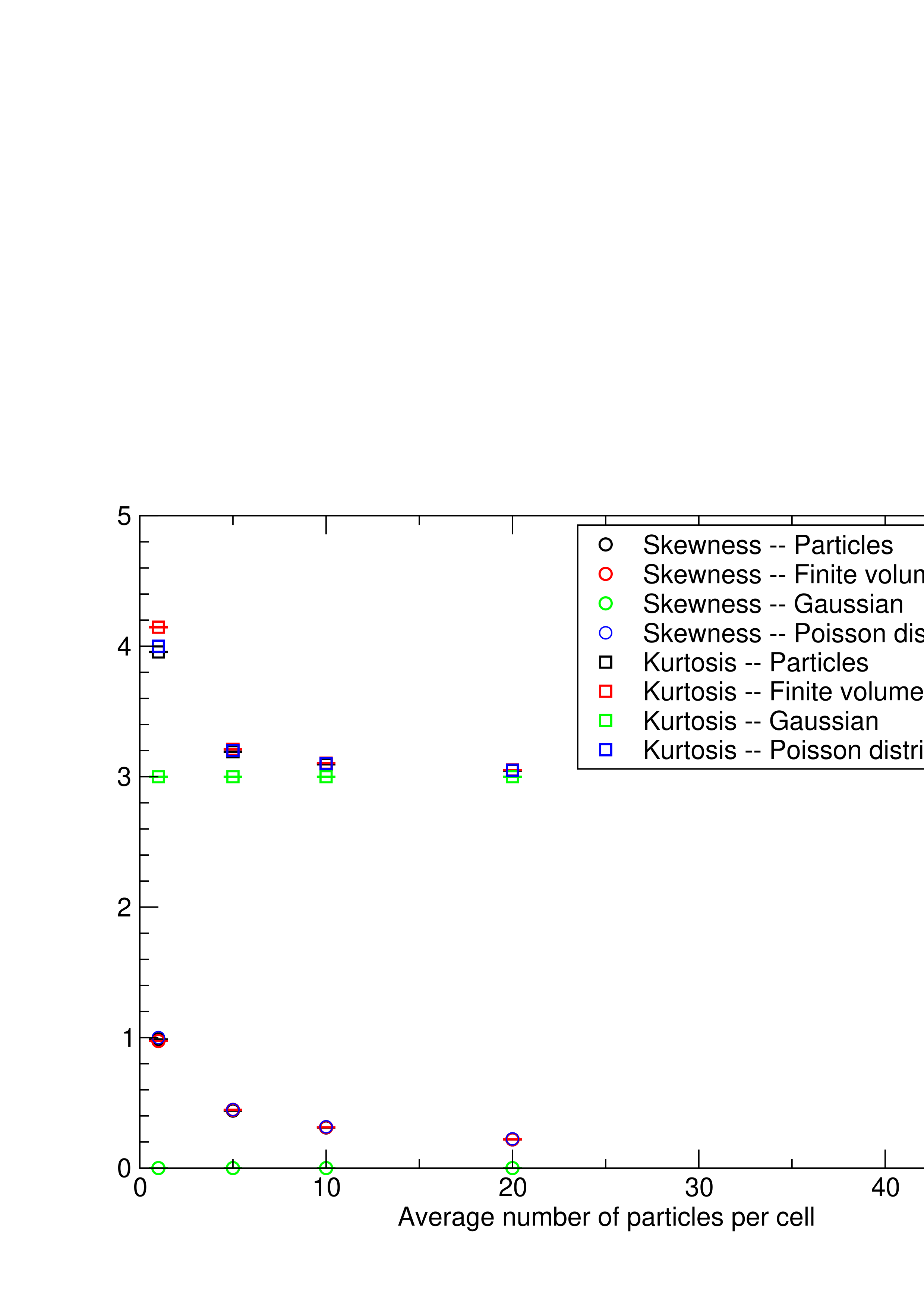}
  \caption{Skewness and kurtosis as a function of particles per cell.  \Add{We have also included skewness and kurtosis for the Poisson distribution.}}
  \label{fig:kurt_skew}
\end{figure}

These higher-order statistics show that accuracy of finite volume SPDE description decreases as the number of particles per \NewOld{cell}{ zone} decreases. This allows us to establish a quantitative relationship between the number of particles per cell and breakdown in the accuracy of the finite volume discretization. 
The data in Figure \ref{fig:kurt_skew} suggest that finite volume  discretization begins to show significant errors with fewer than 5-10 particles per \NewOld{cell}{zone}, suggesting \Add{5-}10 particles per \NewOld{cell}{zone} is a reasonable value to use as a ``refinement" criterion.  \Add{We note that this estimate is somewhat conservative since the criterion is based on the average number of cells rather than the local variability.}
Thus, we expect that the threshold can be set toward the lower end of the range.}

\Add{Refinement criteria based on higher-order statistics is different than what is typically used in hybrid algorithms for fluids  such as those based on DSMC / continuum.  In those cases, criteria based on gradients are typically used, see, for example, \cite{GarciaHornung}.  The use of gradients in the fluid hybrid setting is based on breakdown of the continuum model \cite{Bird_1970,GarciaAlder_1998}. However, tests on the particle system suggest that the SPDE description remains accurate in the presence of strong gradients provided that there are a sufficient number of particles per cell in the region of the gradient.}

\begin{figure}[h!]
  \centering
  a)
  \includegraphics[width=.45\textwidth]{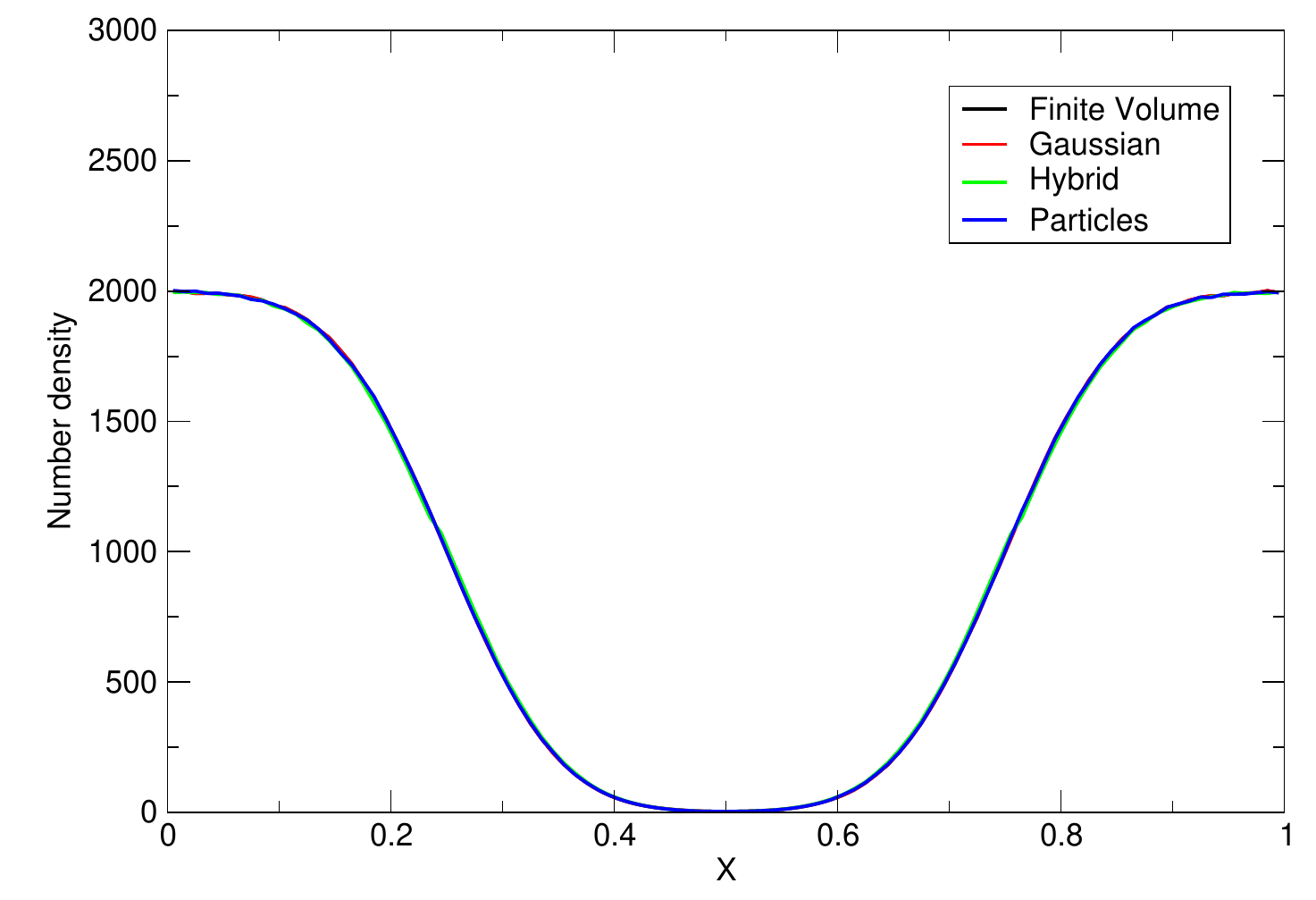}
  b)
  \includegraphics[width=.45\textwidth]{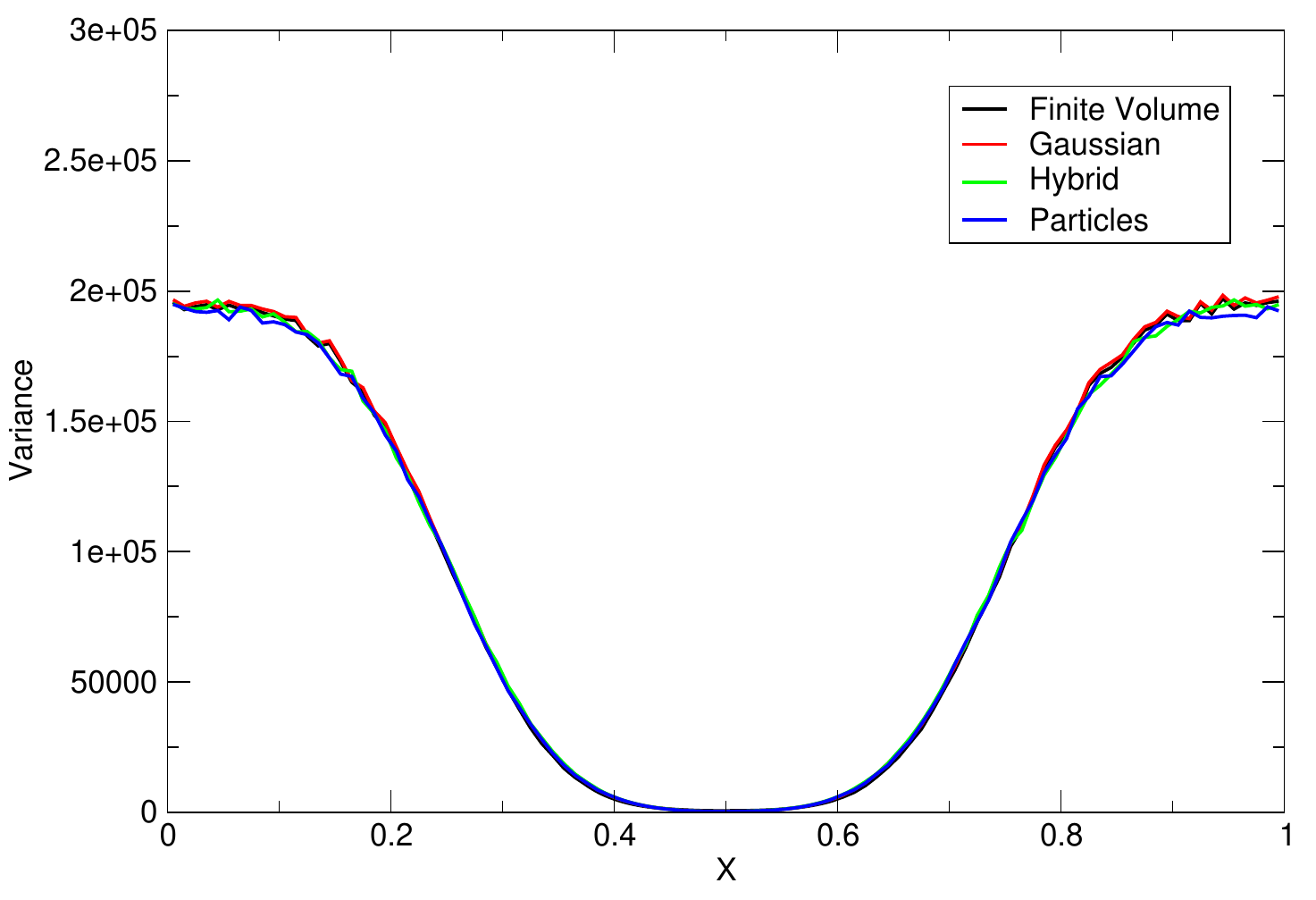}

  c)
  \includegraphics[width=.45\textwidth]{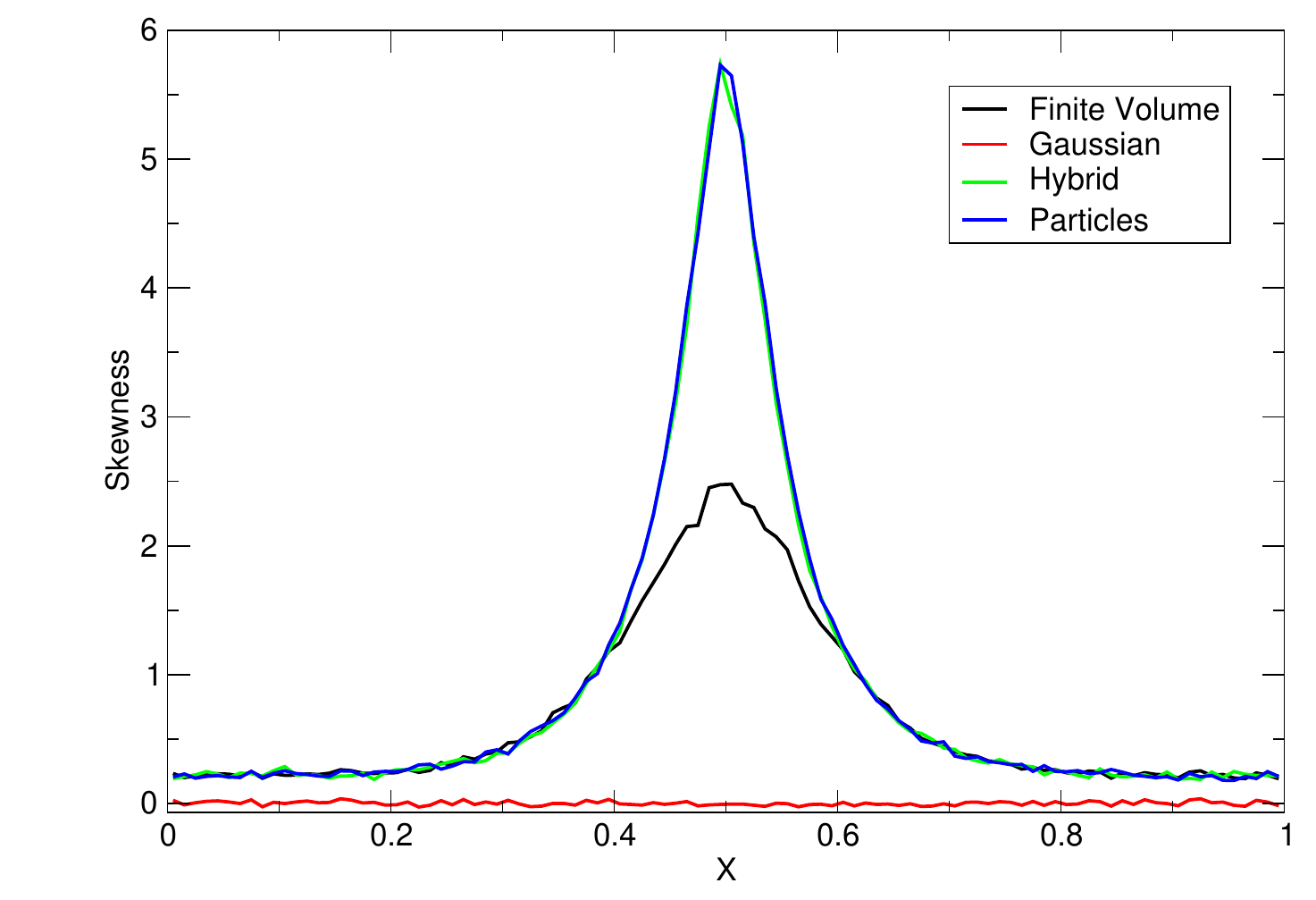}
  d)
  \includegraphics[width=.45\textwidth]{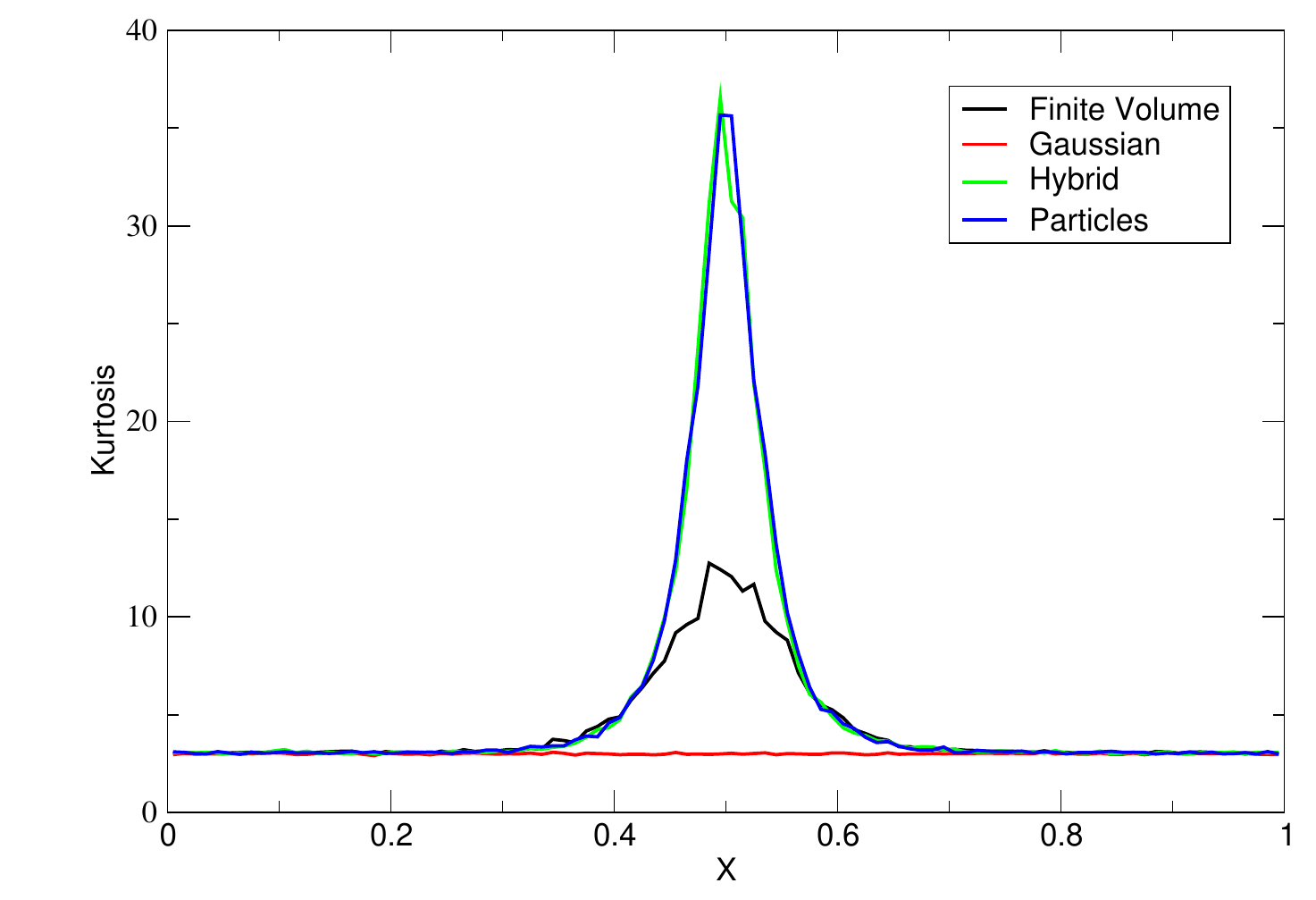}
  \caption{Ensemble average of dynamics of particles diffusing into a void region.  (a) Mean, (b) Variance, (c) Skewness, and (d) Kurtosis.}
  \label{fig:1d_dynamic_stats}
\end{figure}

As an initial test of the hybrid, we consider a dynamic case in which particles are diffusing into an initially void region.
We again consider a one-dimensional periodic case with 100 \NewOld{cells}{ zones} with 
$q = 0$ for $0.25 \leq x \leq 0.75$ and $q=2000$ elsewhere, corresponding to 20 particles per cell.  For this test we fix the particle region to be one-cell wider in each direction than the region where there are initially no particles.
In Figure \ref{fig:1d_dynamic_stats} we present statistics for an ensemble of 20,000 simulations.  The means and variances are in good agreement for all 4 methods.  For the linearized Gaussian method, however, the skewness and kurtosis, which reflect a Gaussian distribution, are significantly different than the values for the particle algorithm, particularly in the low particle density region.  The finite volume algorithm does somewhat better, however, the errors in the low particle density region are still significant. The hybrid algorithm captures the peaks in skewness and kurtosis at low particle densities, demonstrating that the hybrid algorithm can accurately model systems with regions of low particle density.

\begin{figure}[h!]
  \includegraphics[width=.42\textwidth]{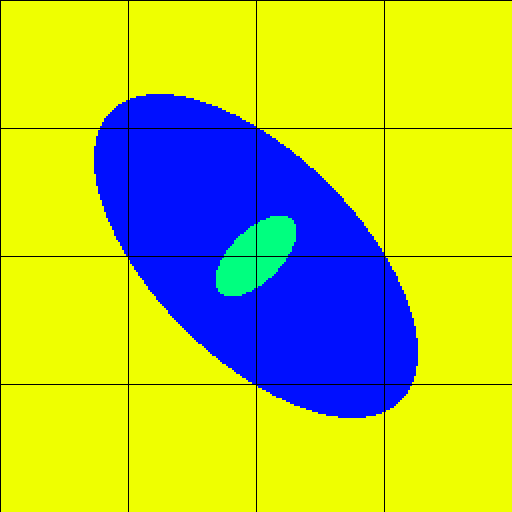}
  \hspace{.05in}
  \includegraphics[width=.42\textwidth]{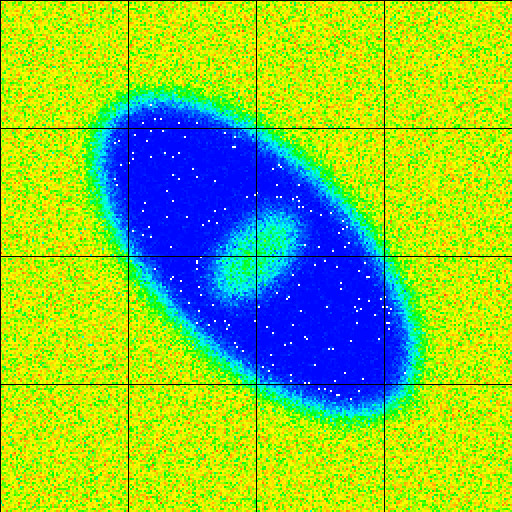}
   \includegraphics[width=.1\textwidth]{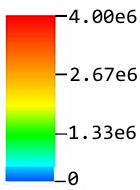}

\hspace{1.2in} $t = 0.$ \hspace{2.3in} $t = 0.00048$

  \vspace{.1in}
 \includegraphics[width=.42\textwidth]{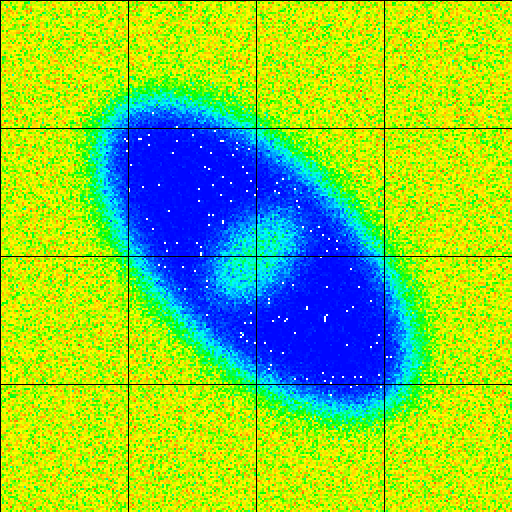}
 \hspace{.05in}
  \includegraphics[width=.42\textwidth]{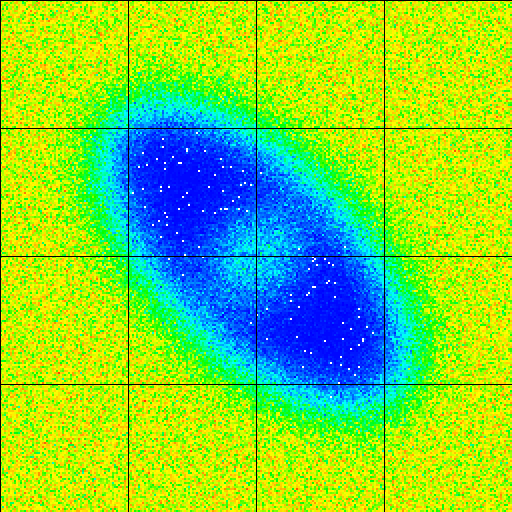}

\hspace{1in} $t = 0.00096$ \hspace{2.2in} $t = 0.0024$
    \caption{\Add{Time sequence for the finite volume method.  White dots  indicate locations where the finite volume solution is negative. The color bar indicates number density.}}
  \label{fig:2d_sequence_fv}
\end{figure}

\begin{figure}[h!]
  \includegraphics[width=.42\textwidth]{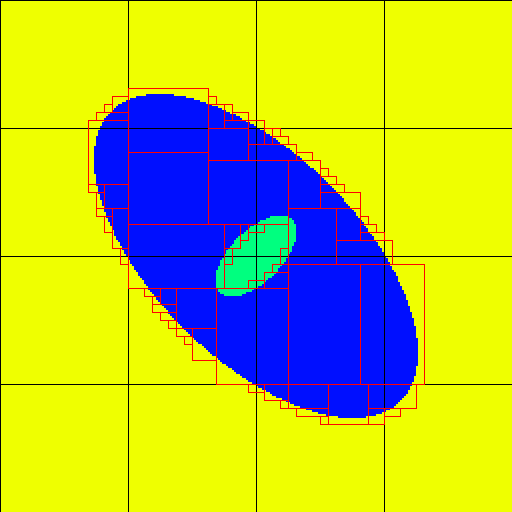}
  \hspace{.05in}
   \includegraphics[width=.42\textwidth]{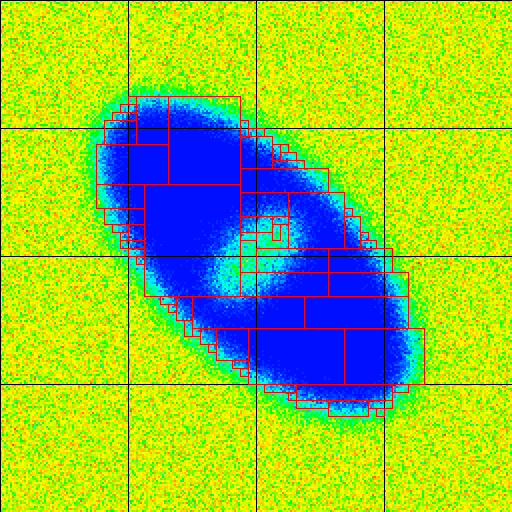}
    \includegraphics[width=.1\textwidth]{colormap_2d.pdf}
    
\hspace{1.2in} $t = 0.$ \hspace{2.3in} $t = 0.00048$

  \vspace{.1in}
   \includegraphics[width=.42\textwidth]{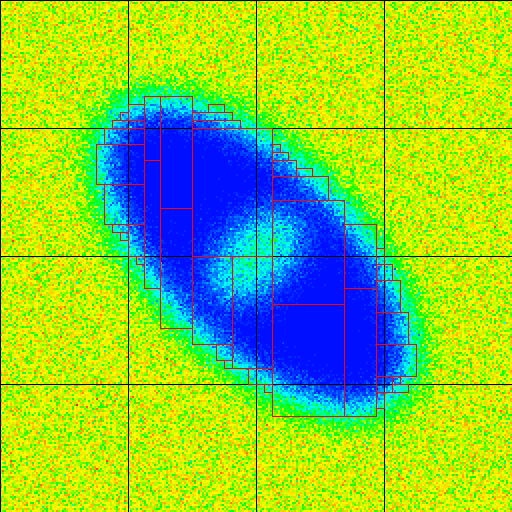}
   \hspace{.05in}
  \includegraphics[width=.42\textwidth]{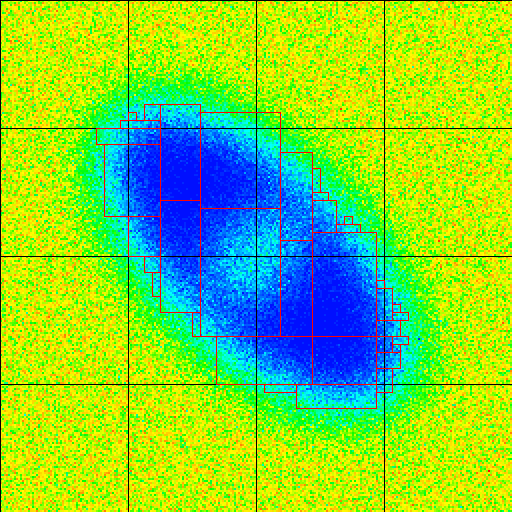}
    
\hspace{1in} $t = 0.00096$ \hspace{2.2in} $t = 0.0024$
  \caption{\Add{Time sequence for the hybrid algorithm.  The red boxes indicate the dynamics location of the particle regions.  The hybrid solution remains positive. The color bar indicates number density.}}
  \label{fig:2d_sequence_hy}
\end{figure}

Next we consider examples in multiple dimension that demonstrate \NewOld{both the dynamic adaptivity and the positivity preservation property of the methodology.}{the dynamic adaptation of the methodology.}  The first example is a two-dimensional case with initial conditions given by an ellipsoidal region with 15 particles per cell surrounded by a larger elliptical region without any particles embedded in a background with 30 particles per cell, as illustrated in Figure \ref{fig:2d_sequence_fv}.  
\Add{
Simulations for the finite volume method and the hybrid algorithm were performed on a $256  \times 256$ grid.  (We note that at this resolution 15 particles per cell corresponds roughly to a number density of $10^6$.)
For the hybrid simulation we tag cells with fewer than 5 particles per cell to be in the particle region, which is near the low end of the range of criteria discussed earlier.
We  present time sequences for finite volume discretization and the hybrid in Figures \ref{fig:2d_sequence_fv} and \ref{fig:2d_sequence_hy}, respectively. 
}
\NewOld{In the hybrid algorithm, the}{The} inner ellipse is initially outside the particle region but as the system evolves the density decreases and the particle region expands to include the entire inner ellipse. At the outer edge of the larger ellipse, particles diffuse into the region and the particle region shrinks.  Thus, this example \NewOld{illustrates the dynamic adaptivity of the hybrid with differ parts of the domain entering and leaving the particles region during the simulation.}{shows regions entering and leaving the particle region. } The hybrid algorithm and the finite volume algorithm produce similar large-scale features; however, the finite volume algorithm generates negative values (indicated by white dots in the figure) whereas the solution from the hybrid algorithm remains positive.

\begin{figure}[h!]
  \includegraphics[width=.43\textwidth]{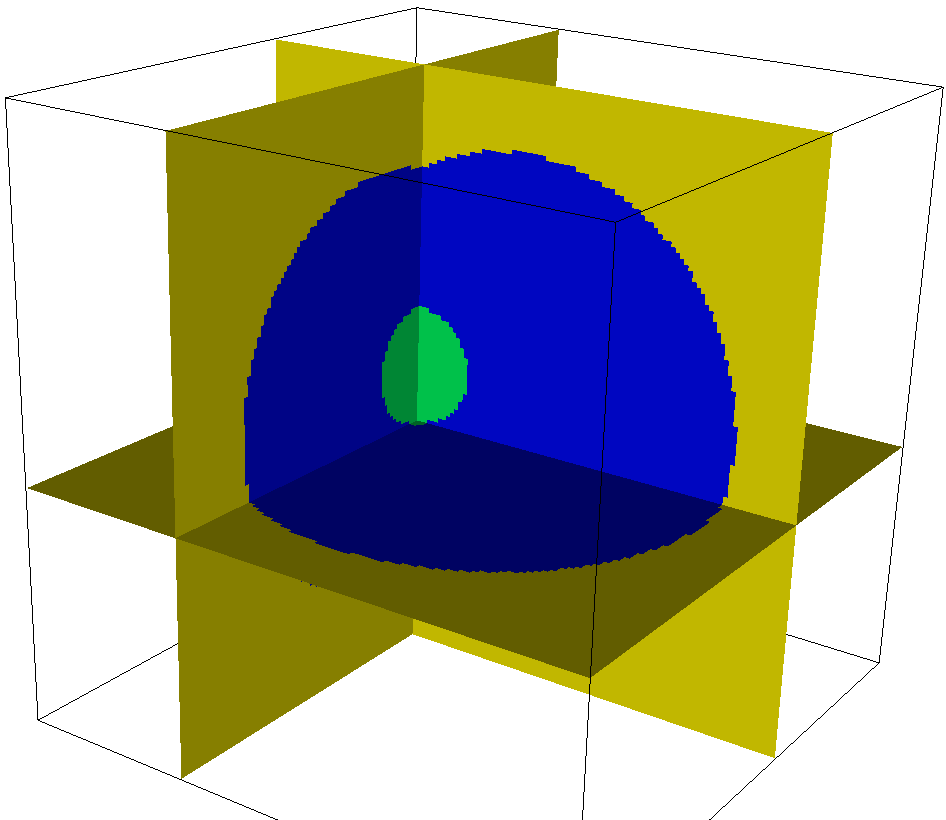}
  \includegraphics[width=.43\textwidth]{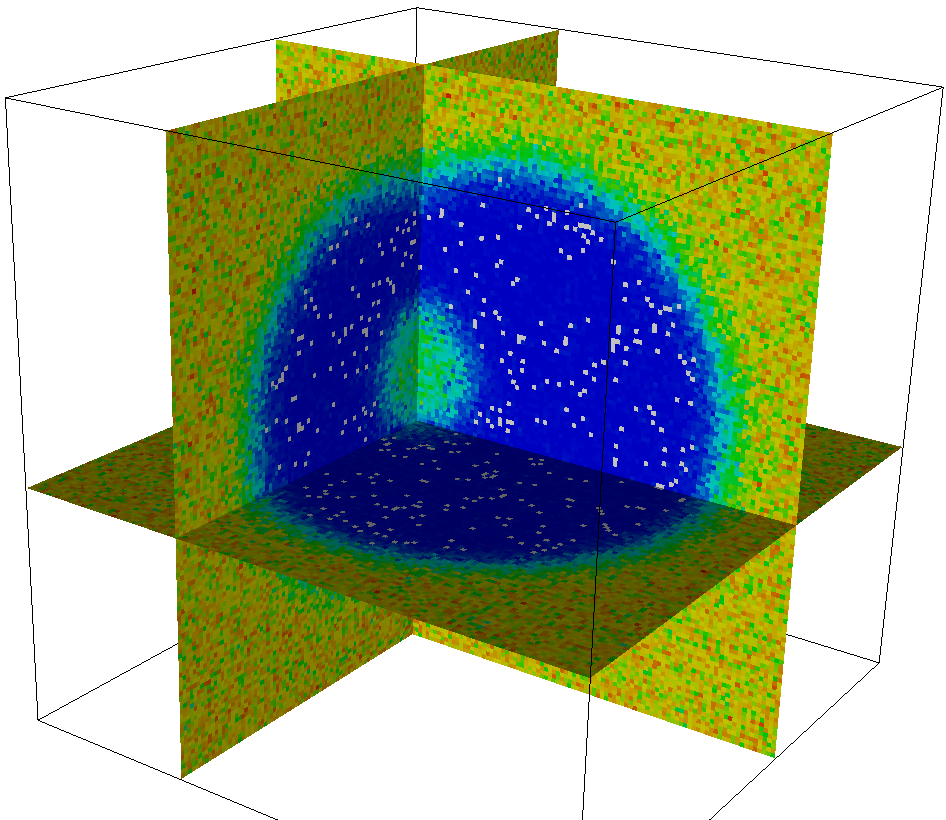}
  \includegraphics[width=.1\textwidth]{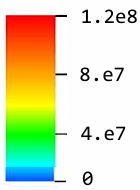}
    
\hspace{1in} $t = 0.0$ \hspace{2.3in} $t = 0.001$

  \vspace{.1in}
    \includegraphics[width=.43\textwidth]{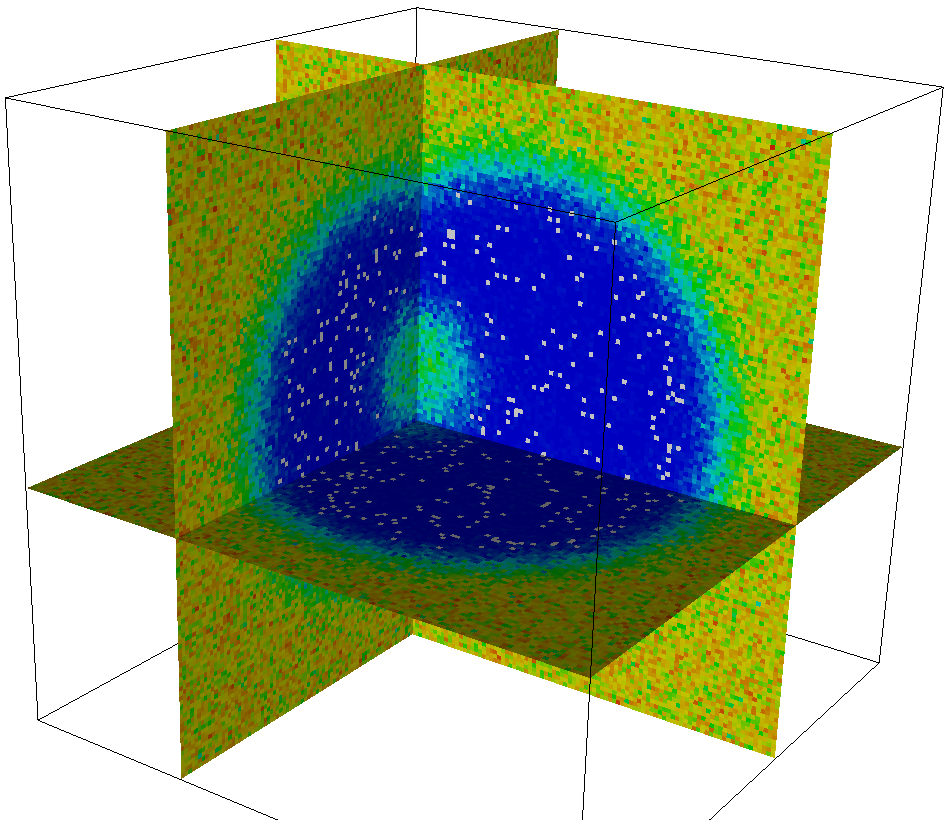}
  \includegraphics[width=.43\textwidth]{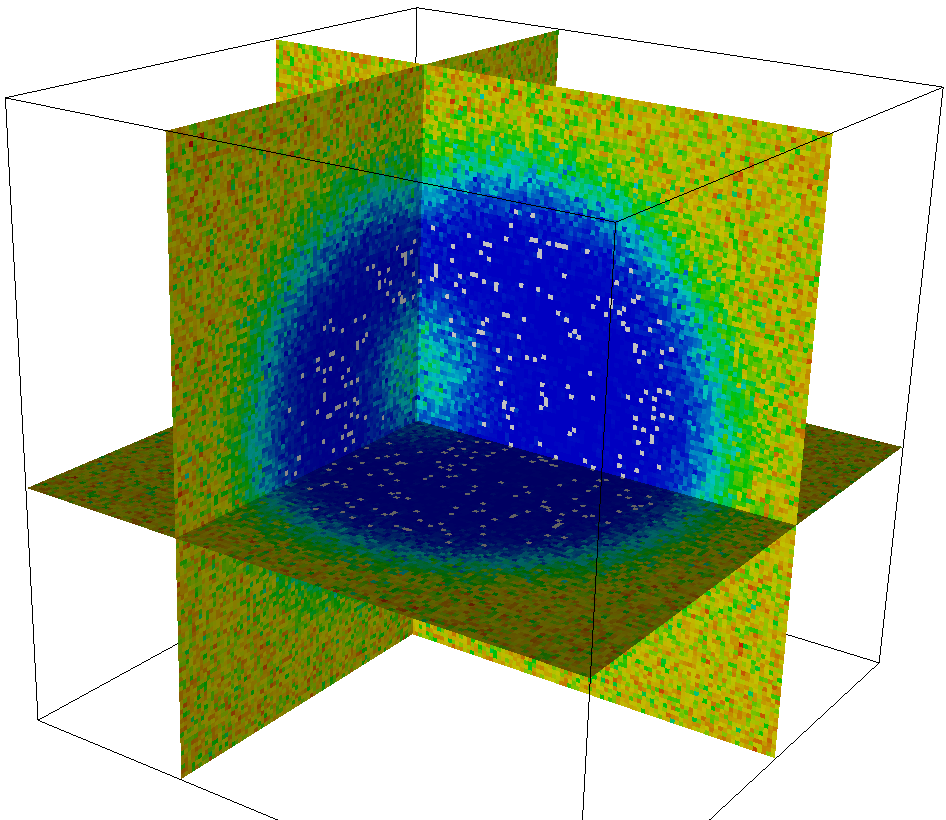}
        
\hspace{1in} $t = 0.0015$ \hspace{2.2in} $t = 0.0025$
     \caption{Time sequence \NewOld{for the finite volume method in three dimensions.}{comparing hybrid algorithm (left) with finite volume algorithm (right).  The boxes outlined in black indicate the location of the particle region. } White dots \RevOut{in finite volume algorithm} indicate locations where the \Add{finite volume} solution is negative.  \Add{The color bar indicates number density.}
     }
  \label{fig:3d_sequence_fv}
\end{figure}

\begin{figure}[h!]
    \includegraphics[width=.43\textwidth]{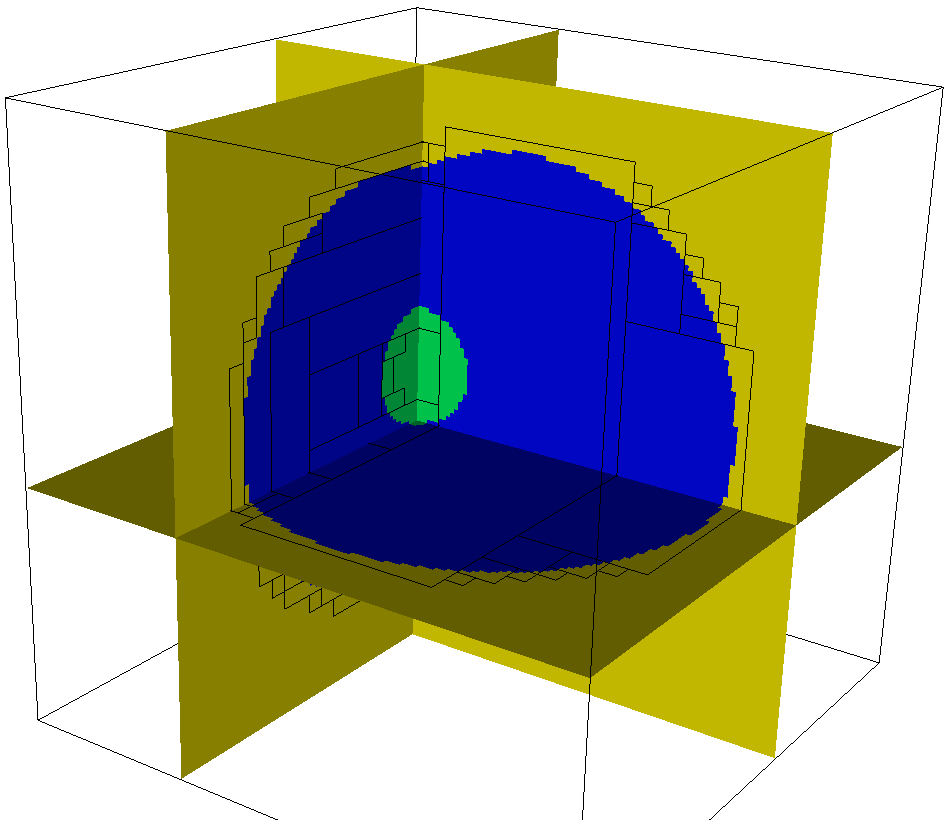}
  \includegraphics[width=.43\textwidth]{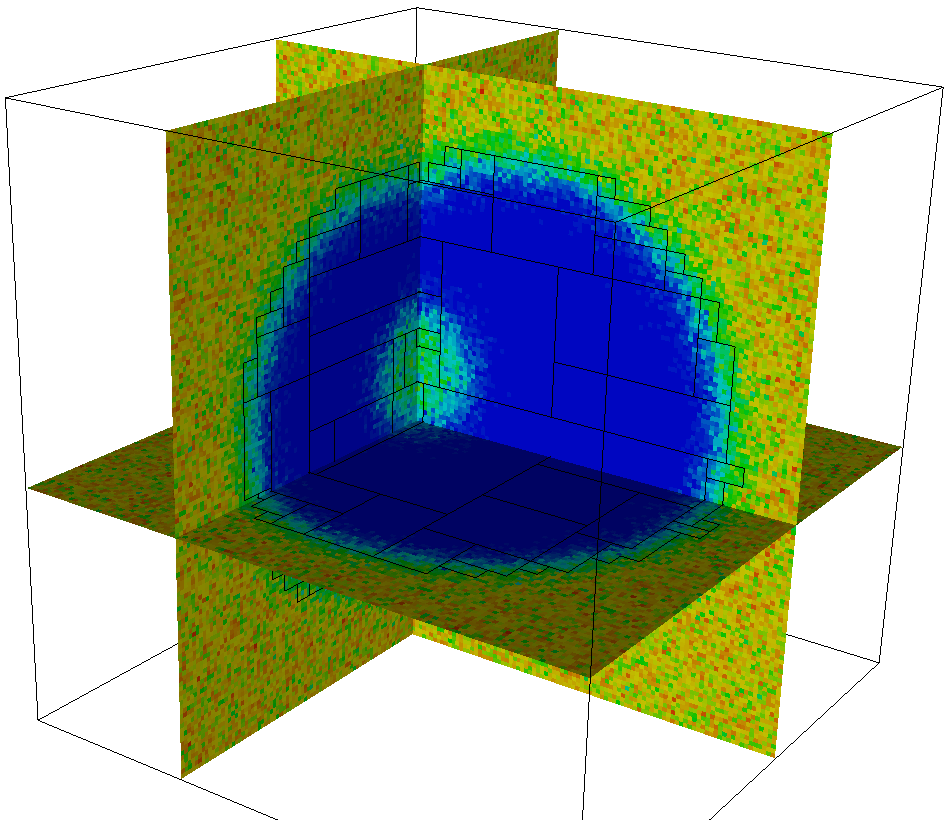}
   \includegraphics[width=.1\textwidth]{colormap_3d.pdf}
    
\hspace{1in} $t = 0.0$ \hspace{2.3in} $t = 0.001$

 \vspace{.1in}
    \includegraphics[width=.43\textwidth]{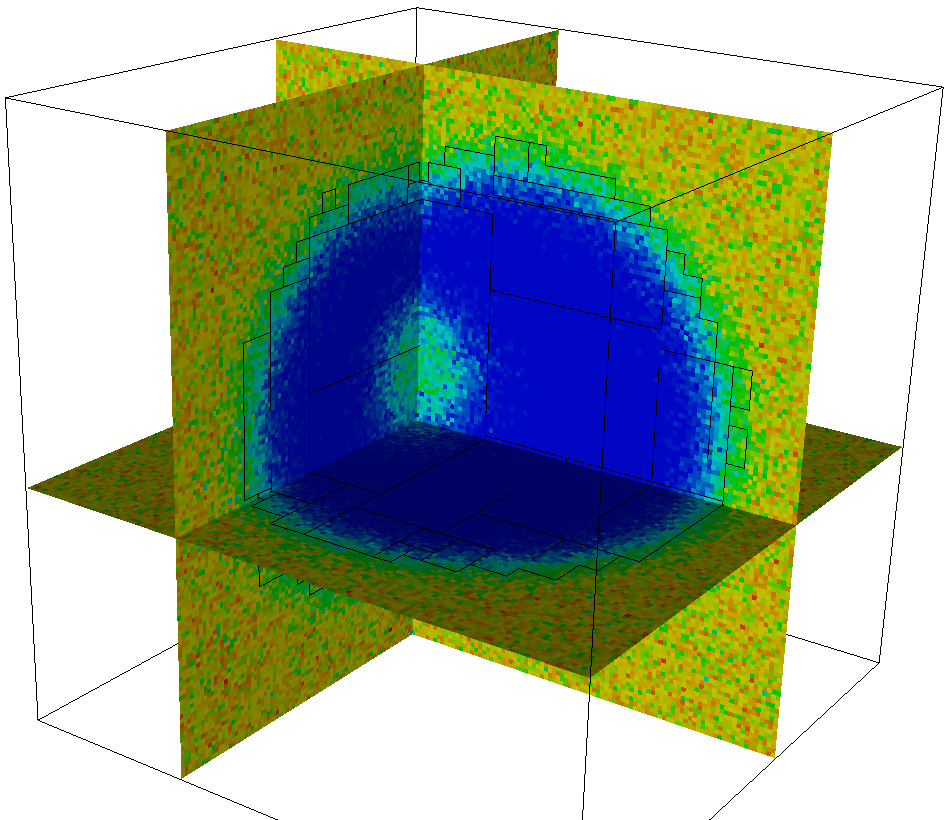}
  \includegraphics[width=.43\textwidth]{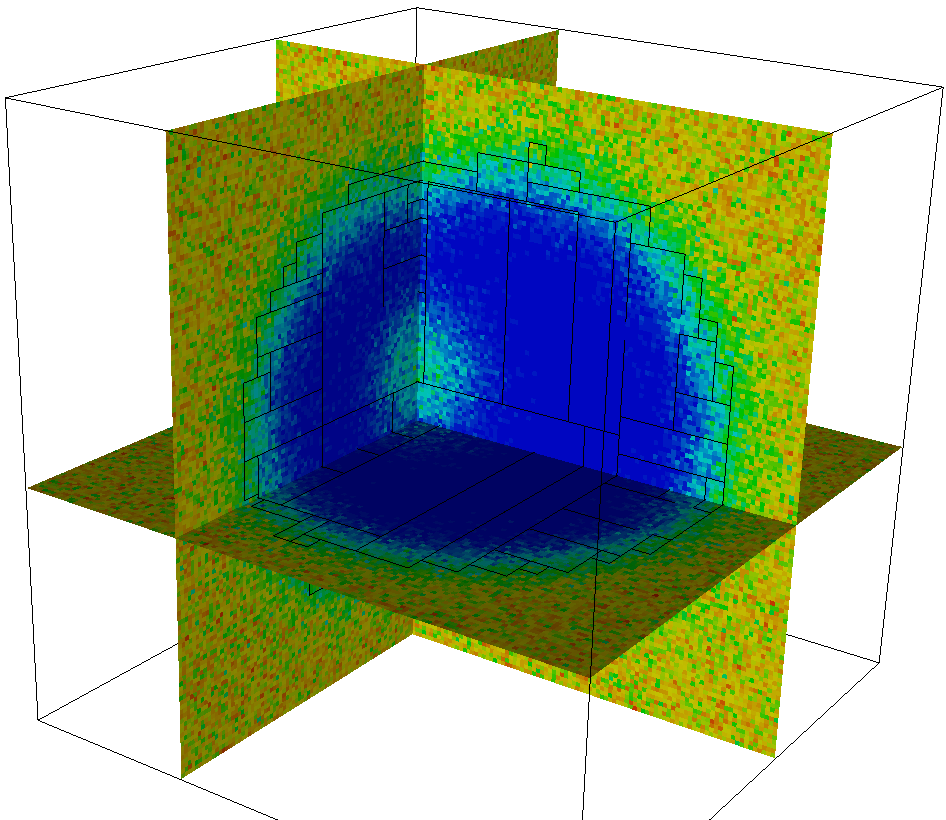}
      
\hspace{1in} $t = 0.0015$ \hspace{2.2in} $t = 0.0025$
 \caption{\NewOld{Time sequence of the hybrid algorithm in three dimensions. The boxes outlined in black indicate the location of the particle region.  The hybrid solution remains positive throughout the simulation. The color bar indicates number density.}{Continuation of time sequence comparing hybrid algorithm (left) with finite volume algorithm (right). The hybrid solution remains positive. Times correspond to $t= 0.0015$, and $t= 0.0025$.}}
  \label{fig:3d_sequence_hy}
\end{figure}

\NewOld{Our next}{The final} example is a
three-dimensional case with initial conditions given by a spherical region with 15 particles per cell surrounded by a larger spherical region without any particles embedded in a background with 30 particles per cell, as illustrated in Figure \ref{fig:3d_sequence_fv}. 
\Add{ 
Simulations with the finite volume method and the hybrid algorithn were performed on a $128 \times 128 \times 128$ grid.  At this resolution 15 particles per cell corresponds to a number density of roughly $3\times 10^7$. As in the two dimensional case, for the hybrid we \Add{again} tag cells with fewer than 5 particles per cell to be in the particle region. Results for the finite volume discretization and the hybrid are presented in Figures \ref{fig:3d_sequence_fv} and \ref{fig:3d_sequence_hy}. 
}
The inner sphere is initially outside the particle region. As the system evolves the density decreases and the particle region begins to expand.  At $t=0.001$ only a small region of the inner sphere is still be modeled with the finite volume method. By $t=0.0015$ the particle region has grown to include the entire inner sphere.  At the outer edge of the larger sphere, particles diffuse into the region and the particle region shrinks. The hybrid algorithm and the finite volume algorithm again produce similar large-scale features; however, the finite volume algorithm generates negative values (indicated by white dots in the figure) whereas the solution from the hybrid algorithm remains positive.

\Add{Our final numerical example illustrates the application of the methodology to a problem with an external potential.  
The modifications required to incorporate an external potential into the algorithm are discussed in the Appendix A.
Here we consider a potential of the form
\[
V(x,y) = (x-\alpha)^2 (x-\beta)^2 + (y-0.5)^4 \;\;\;\; \mathrm{on} \;\;\; [0,1]\times [0,1] \;\;\; ,
\]
which has minima at $A = (\alpha,0.5)$ and $B = (\beta,0.5)$.  When $\gamma$ in Eq. (\ref{eq:langevin}) is sufficiently small, the transition from $\alpha$ to $\beta$ becomes a rare event. 
Problems of this type arise in the calculation of rate coefficients in physics, chemistry and biology, originating with early work of Kramers \cite{KRAMERS1940284}. 
Schuss \cite{schuss1980singular} presents a diffusion model for chemical reactions, where the confining potential well  consists of a series of holes and barriers and the particle dynamics is given by the Langevin dynamics of the type
\eqref{eq:langevin}. 
H\"anggi {\it et al.} \cite{Hanggi_1990} discusses a number of developments in this area.  
A seminal paper of Jordan, Kinderlehrer and Otto  \cite{jordan1998variational} investigates Fokker–Planck equations where the drift term is given by the gradient of a potential, which corresponds to the dynamics described by \eqref{eq:langevin}. 
Similar types of problems also arise in large-deviation theory, see Berglund \cite{berglund2011kramers} and are directly relevant to simulations in molecular dynamics \cite{schutte2023overcoming}.

We set $\alpha = 0.3$ and $\beta = 0.7$ and consider the case where $\gamma = 5.0 \times 10^{-4}$.
We initialize the simulation with roughly $10^5$ particles in the subregion $[0.2,0.4] \times [0.25,0.75]$ distributed based on the equilibrium distribution
\[
\frac{1}{\mathcal{Z}} e^{-2 V(x,y) / \gamma}
\]
restricted to the subregion ($\mathcal{Z}$ is the normalization constant).
We discretize the system on a $100 \times 100 $ grid with homogeneous Dirichlet boundary conditions and simulate for 100,000 time steps with $\Delta t = 3.90625 \times 10^{-7}$ using the finite volume scheme, particles dynamics and the hybrid method.  See Appendix B for the treatment of boundary conditions.  Results from the three methods are shown in Figure \ref{fig:exp}.  The presence of a potential leads the bulk of the particles being localized to small regions of the
domain with very low particle density in the remainder of the domain. This results in a loss of positivity for the finite volume scheme in a significant fraction of the domain as indicated in Figure \ref{fig:exp}.

 \begin{figure}[h!]
   \includegraphics[width=.28\textwidth]{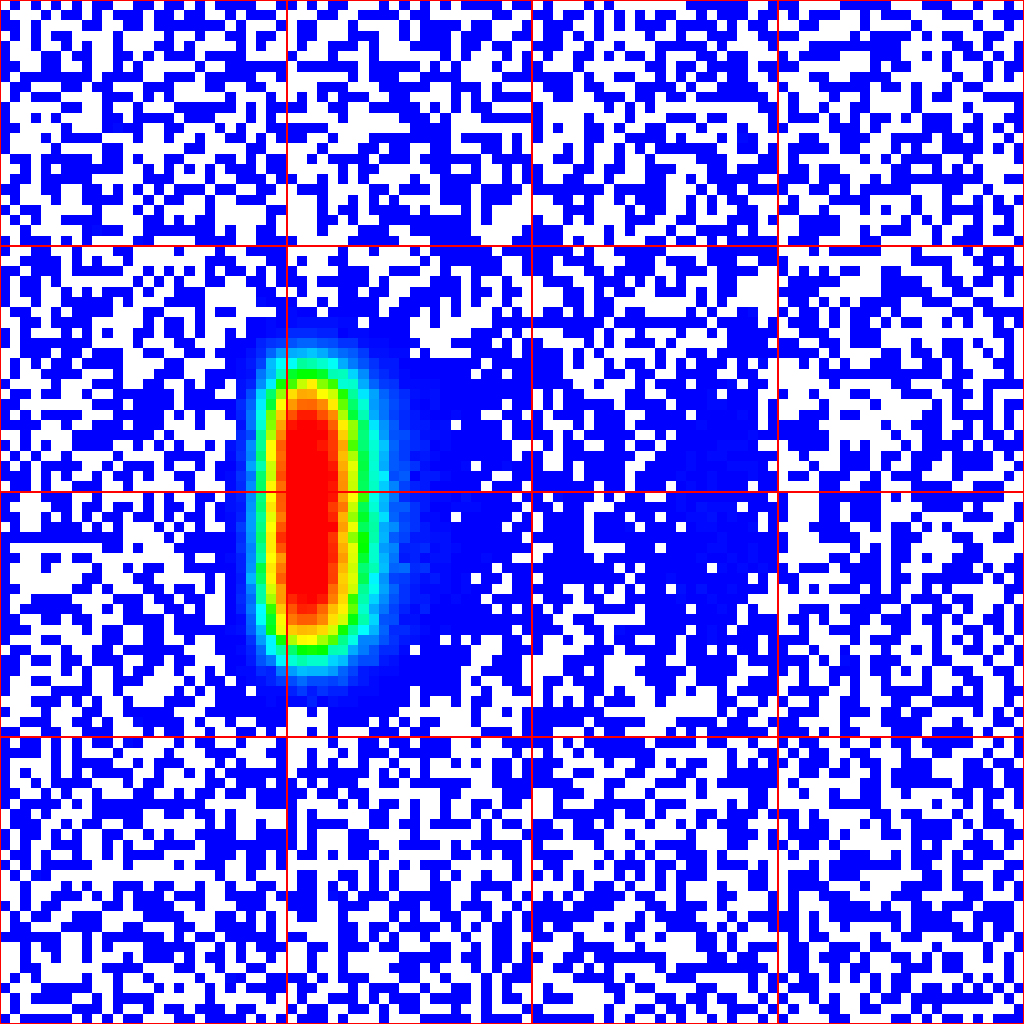}\hspace{.05in}
  \includegraphics[width=.28\textwidth]{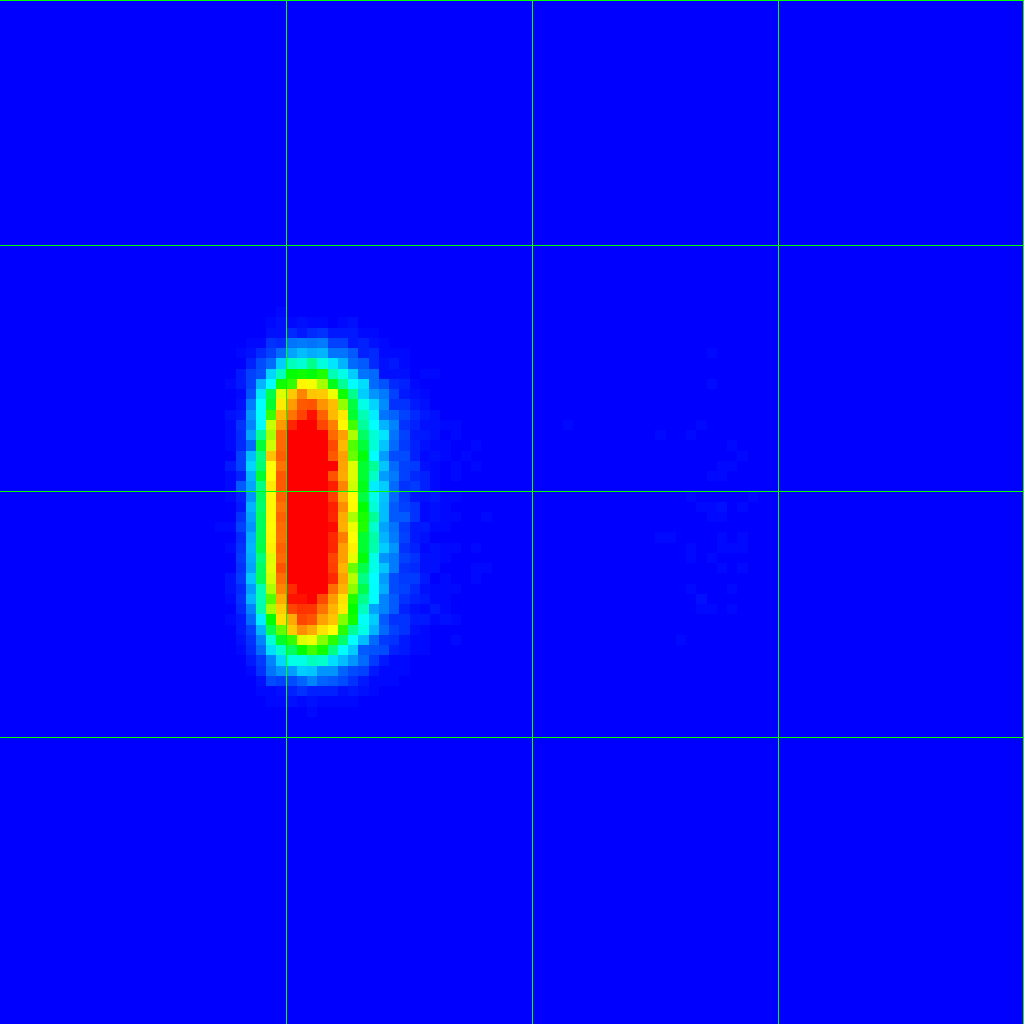}\hspace{.05in}
   \includegraphics[width=.28\textwidth]{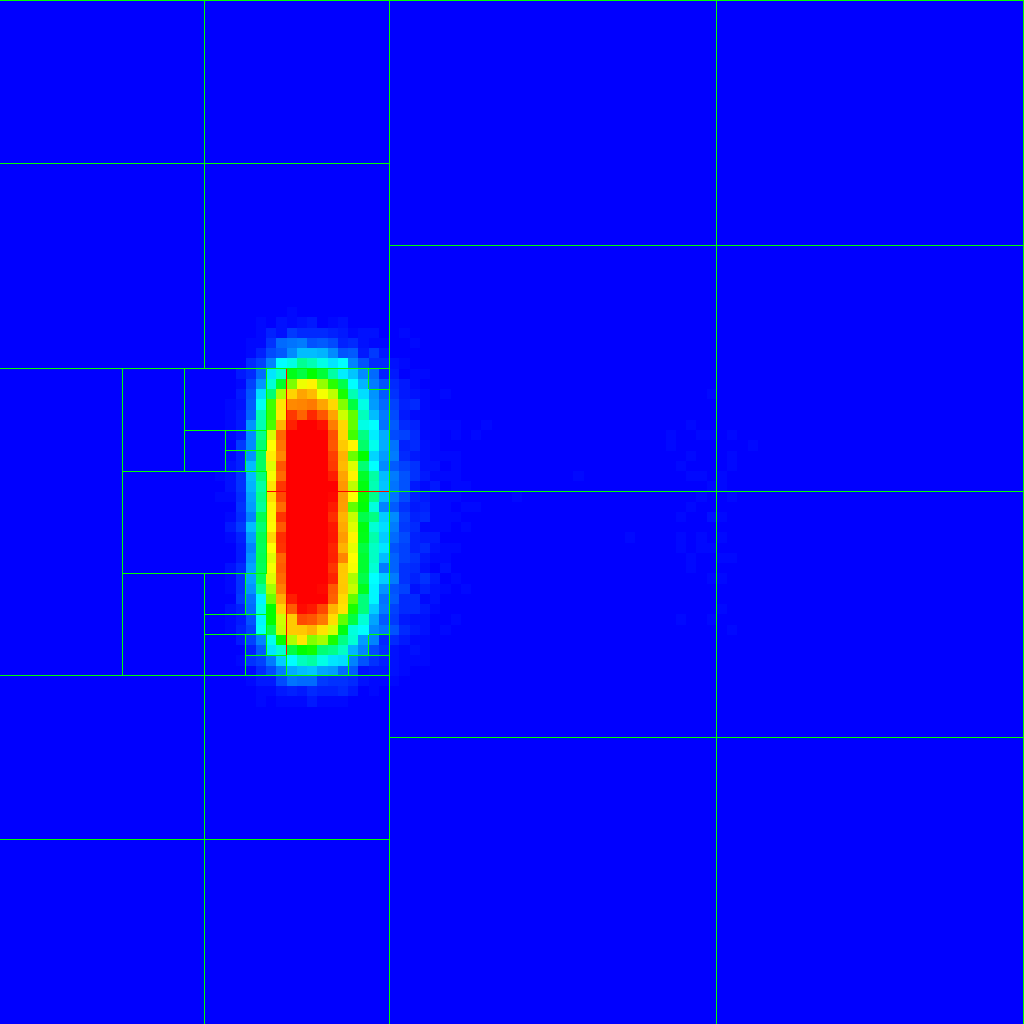}
    \includegraphics[width=.1\textwidth]{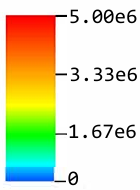}

  \caption{\Add{Comparison of simulation with an external potential at $t = 3.90625\times 10^{-2}$. The left frame is the finite volume algorithm, the center frame is the particle algorithm and the right frame is the hybrid.  White points indicate where the solution has become negative in the finite volume solution.  The light blue lines in the hybrid plot show the grids that define the particle region.  Although most of the domain is in the particle region, there are very few particles in the particle region.}}
  \label{fig:exp}
\end{figure}

In Figure \ref{fig:exp_zr}, we reduce the range on the images to highlight the particles that have crossed the energy barrier.  We consider particles whose $x$ location is greater than 0.6 to be in the potential well at $B$.  Here 0.6 corresponds to the location  where the potential is half the barrier height.  We measured the number of particles that had crossed the barrier into the potential well at $B$ for a small ensemble of ten simulations with each method.  For the particle algorithm and the hybrid method we obtained means at the final time of 388 and 384 particles respectively, which are within statistical error bars for the small ensemble.  The finite volume algorithm gave a mean of 370 particles, which is a statistically significant under prediction.  A much larger ensemble from a one-dimensional case showed similar results with the finite volume method under predicting the number of particles crossing the barrier while hybrid and particle simulations predict means that differ by less than 1\%.

As noted above, when we generate particles from the SPDE representation during regridding or when filling ghost cells, we do not assume that the particles are uniformly distributed in the cell when there is an external potential. Instead we sample grid locations from a local approximation to the equilibrium distribution.
To illustrate the importance of using the appropriate measure for this operation, we performed a small ensemble of simulations with the external potential in which
we defined particle locations for regridding and filling ghost cells using a uniform distribution.  For this ensemble, the mean number of particles that crossed to 
the potential well at $B$ was 555, representing a dramatic over-estimation of the number of particles crossing the potential barrier.    

 \begin{figure}[h!]

     \includegraphics[width=.28\textwidth]{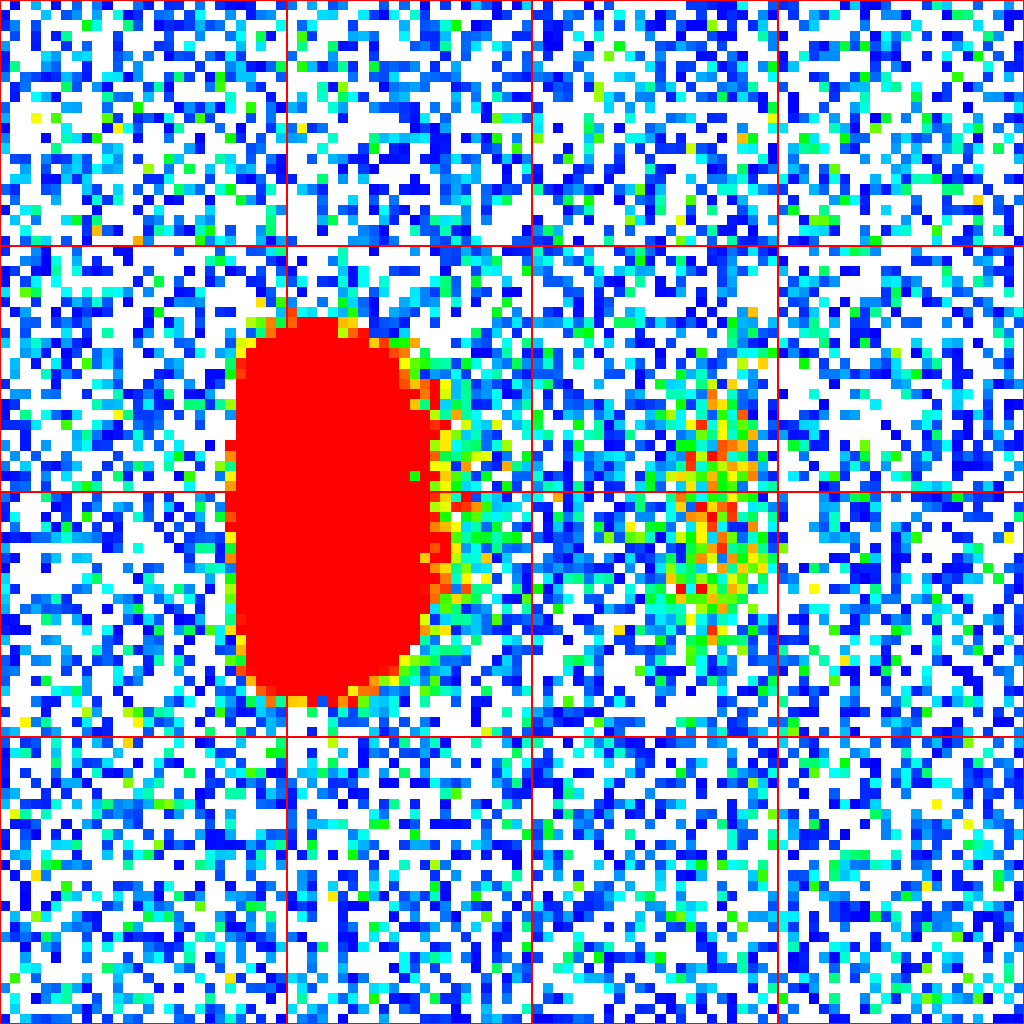}\hspace{.05in}
  \includegraphics[width=.28\textwidth]{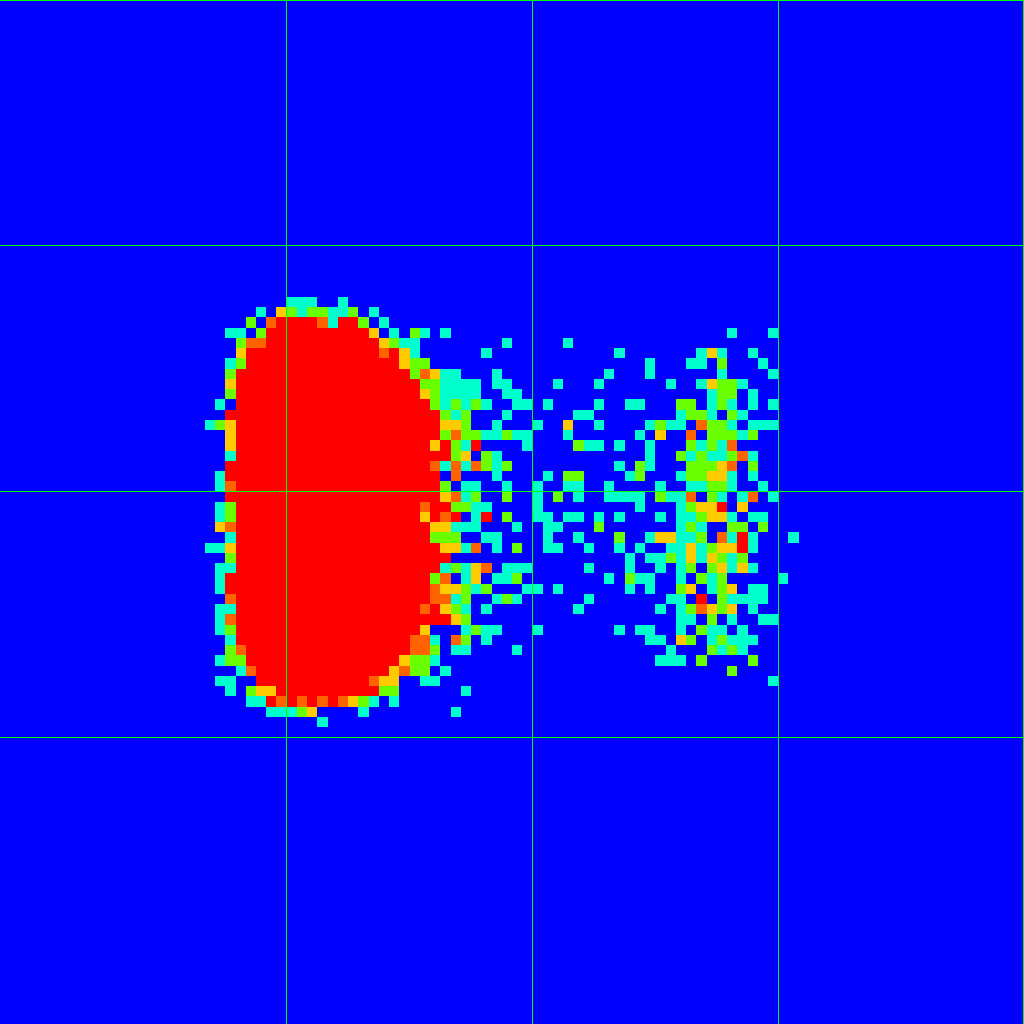}\hspace{.05in}
   \includegraphics[width=.28\textwidth]{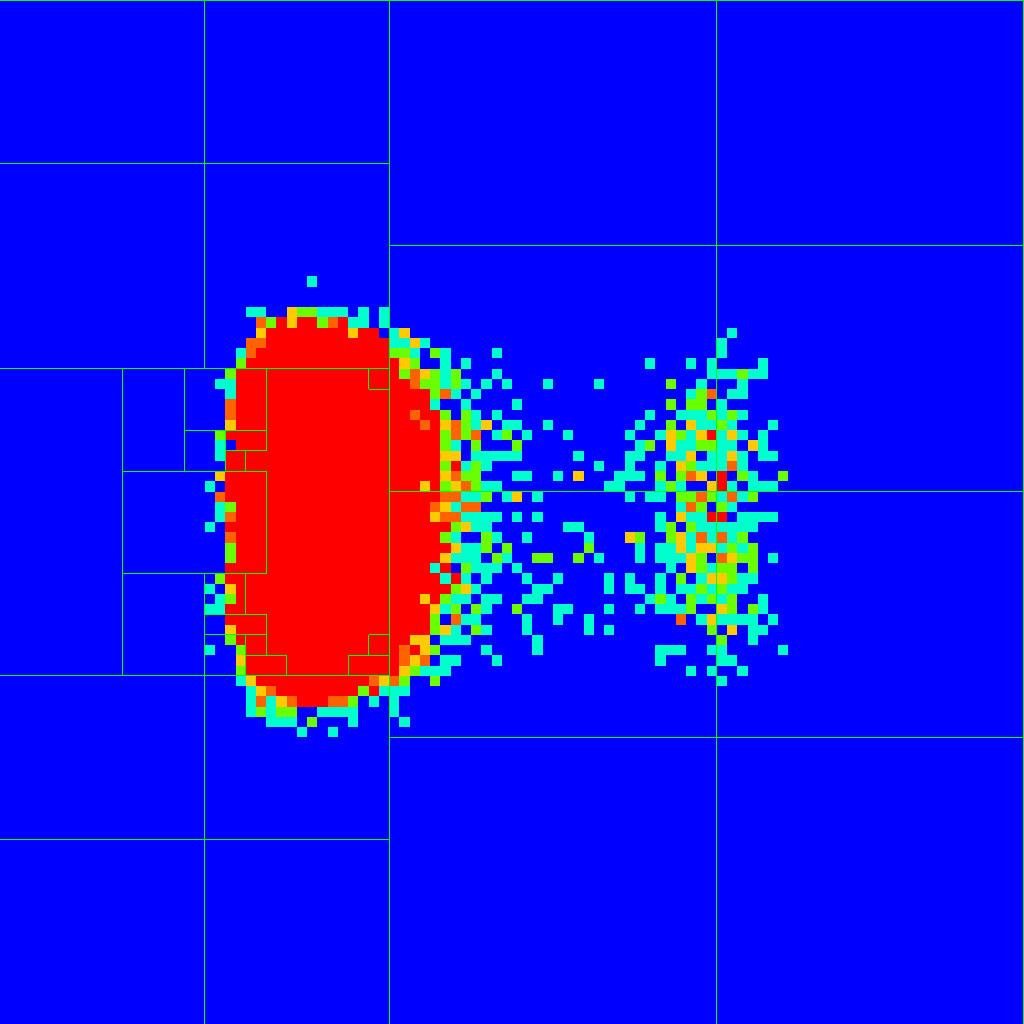}
    \includegraphics[width=.1\textwidth]{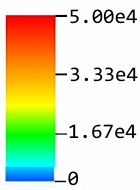}
   
  \caption{\Add{External potential simulations with color ranged reduced by a factor of 100. As above, the left frame is the finite volume algorithm, the center frame is the particle algorithm and the right frame is the hybrid. }}
  \label{fig:exp_zr}
\end{figure}

}

\section{Conclusions and further work}

In this work we have presented a hybrid method for simulating the dynamics of non-interacting particles \Add{with the example that also includes external potential}. This approach combines the simulation of the Dean-Kawasaki equation with truncated noise and particle dynamics in low-density regions. In particular, we developed the refinement criteria by comparing higher-order statistics. Compared to prior work on hybrid algorithms for particle systems, such as \cite{ALEXANDER2002,alexander2003algorithm}, we focused specifically  on regimes where particle density becomes low. In these regimes, our numerical experiments have shown that incorporating a particle description of the system is crucial for accurately preserving higher-order statistics and positivity of the solution. Additionally, we have extended the method to higher spatial dimensions with dynamic adaptation \Add{and we have provided a numerical example that includes an external potential}.

\Add{
An important future direction for this work is the extension of the hybrid framework to include interactions between the particles.  
There has been a substantive amount of work on SPDE models that include interactions between the particles, see, for example, \cite{wehlitz2024approximating, helfmann2021interacting,illien2024dean,bressloff2024generalized}; however, there has been essentially no work on development of hybrid algorithms of the type considered here.
Developing hybrids that include particle interactions would enable a broad range of applications in physics, chemistry, biology and
the social sciences. 
Unfortunately, incorporating these types of models into a hybrid poses a number of challenges.  First, one needs to be able to compute the interaction between parts of the solution represented by particles and parts of the solution represented by number density.  The second major issue is how to generate particles consistent with the SPDE state for regriding and to fill ghost cells.  As shown here for an external potential, the numerical results for the hybrid are extremely sensitive to how particles are generative. 
The solution to either of these problems is not known.

Another promising application of the methodology would be to study particles immersed within a surface  for  problems where geometry is crucial, as discussed in \cite{rower2022surface}. Here the goal would be to develop a hybrid method that incorporates the geometry by simulating the dynamics of Brownian particles on a surface, solving the corresponding stochastic equation on that surface, and synchronizing these processes while accounting for the geometric factors.

There are also a number of more mathematical research topics suggested by this work.
One area would be to prove error bounds for the FVM applied to the modified equation from \cite{djurdjevac2022weak} and obtain a rate of convergence. Given that the modified equation has a well-behaved solution, one could directly analyze the error between the moments of the discretized and continuous solution, rather than using a weaker norm. The goal would be to compare how the choice of parameters determined by the discretization error aligns with the parameters chosen in \cite{djurdjevac2022weak} for the error between the modified and Dean-Kawasaki equations, as well as with the parameters obtained from simulations using the hybrid method. This type of analysis would allow us to link error estimates for the hybrid to the refinement criteria.

There are also opportunities to improve the basic SPDE methodology.  Of particular interest would be higher-order time discretization methods, such as the Milstein scheme or the BDF2 method. Note that BDF2 is  expected to provide higher accuracy for SPDEs with small noise, as discussed in \cite{buckwar2006multistep}. Since the noise factor in Eq. \eqref{eq:DK_with_reg_noise} is $N^{-1/2}$, we expect to achieve higher-order accuracy when the number of particles is sufficiently large. It would be interesting to compare this threshold with the number of particles per cell required to ensure correct higher-order statistics and investigate the whether improvements to the temporal integrator impacts the refinement criterion.

}

\Add{
\section*{Appendix A:  Discretization including an external potential}

In this appendix we briefly discuss the modifications needed to include an external potential into the dynamics.
We will retain the simple Euler-Maruyama temporal discretization that was used in the absence of an external potential. The discretization of the particle motion then becomes
\[
X_i^{n+1} = X_i^n -\frac{\Delta t}{\gamma} \nabla V(X_i^n) + \sqrt{\Delta t} N_i^d
\]
where $N_i^d$'s are vectors of independent sample from a normal distribution.

For the SPDE, we add an additional flux corresponding to the potential term $\nabla \cdot (u \nabla V)$.  Specifically, we define
\[
F^{pot,x}_{i+\half,j,k} = \frac{q_{i,j,k}^n+q_{i+1,j,k}^n}{2 \gamma}\; \nabla V(x_{i+\half},y_j,z_k).
\]
We then define the total flux to be
\[
F_{i+\half,j,k}^x = \detF_{i+1/2,j,k}^x + \stochF_{i+1/2,j,k}^x + F^{pot,x}_{i+\half,j,k} .
\]

With these modifications, the hybrid algorithm is essentially unchanged.  The only modification that is needed is in the algorithm for defining a particle presentation in a cell that is consistent with the SPDE solution. 
When there is an external potential, the equilibrium distribution of particles is no longer uniform.  Instead it is given by
\begin{equation}
    p(X) = \frac{1}{\mathcal{Z}} e^{-2 V(X)/ \gamma}
    \label{eq:eq_dist}
\end{equation}
where $X$ is the position in $d$ space and $\mathcal{Z}$ is a normalization factor. (For the conditions considered here, $V$ is large enough at the domain boundary that the support of $p$ is essentially contained inside the domain.)
When we need to define particles in a given cell consistent with the SPDE solution in that cell, we want to conditionally sample form Eq. (\ref{eq:eq_dist}) instead of assuming the particles are uniformly distributed.  This can be important when the potential changes significantly within a given cell.  Here, we have constructed a local approximation to equilibrium potential based on the potential gradient at the cell center. This approximation allows us to account for variability within the cell but is straightforward to sample.

For both particle update and SPDE integration, one needs to consider whether that addition of the external potential necessitates a smaller time step. In particular, we still need to constrain a particle to not move too far during a time step and ensure the time is small enough that the SPDE discretization is stable and accurate.  For the case considered here we have not found  reduction in time step  necessary.

\section*{Appendix B:  Discretization of boundary conditions}

For Dirichlet boundary conditions, for the SPDE we impose the boundary condition on the face that coincides with the boundary.  The diffusive component of the flux is then
given by 
\[
\overline{F}_{B,j,k}^x = \frac{q_B^n - q_{i,j,k}^n}{\Delta x} \;\;\;\mathrm{etc.},
\]
where $q_B^n$ is the Dirichlet value on the boundary at $t^n$. 
This formula reflects the reduction in the distance between the points to $\Delta x / 2$.
The stochastic flux at the boundary is then multiply by $\sqrt{2}$ to preserve fluctuation balance at the boundary.
The potential flux is simply evaluated using the boundary condition.

For the particle regions that touch the boundary, the boundary condition is represented by a reservoir of particles outside the domain, as done in \cite{ALEXANDER2002}. The number of particles in the reservoir is obtained by sampling from a Poisson distribution with mean given by the boundary condition. The location of the particles in the reservoir are sampled from the appropriate equilibrium distribution.  After the particles are advanced, particles outside the domain are discarded.

}

\section*{Acknowledgement}
 ADj gratefully acknowledges funding by Daimler and Benz Foundation as part of the scholarship program for junior professors and postdoctoral researchers. 
 Further support of ADj is by Deutsche Forschungsgemeinschaft (DFG)  through grant CRC 1114 "Scaling Cascades in Complex Systems", Project Number 235221301, Project C10 "Numerical
Analysis for nonlinear SPDE models of particle systems" and by a Hanna Neumann Fellowship of the Berlin Mathematics Research Center Math+.
The work of AA and JB was supported by the U.S. Department of Energy, Office of Science, Office of Advanced Scientific Computing Research, Applied Mathematics Program under contract No. DE-AC02-05CH11231.

\bibliography{main_2}
\bibliographystyle{abbrv}

\end{document}